\documentclass[aos,MSNbibl,seceqn,citesort,rotating,dvips]{arximspdf}
\usepackage{graphicx}

\doi{10.1214/11-AOS908}
\volume{39}
\issue{5}
\pubyear{2011}
\firstpage{2502}
\lastpage{2532}

\makeatletter

\newproclaim{remark}{Remark}
\newtheorem{theorem}{Theorem}[section]

\newcommand{\uU}{\tilde{U}}

\newcommand{\be}{\beta}
\newcommand{\Bigmi}{\mid}
\newcommand{\bU}{{\bar U}}
\newcommand{\bze}{{\bar\ze}}
\newcommand{\cF}{\mathcal{F}}
\newcommand{\De}{\Delta}
\newcommand{\ep}{\varepsilon}
\newcommand{\ga}{\gamma}
\newcommand{\half}{^{1/2}}
\newcommand{\ib}{_{i\bullet}}
\newcommand{\ourij}{_{ij}}
\newcommand{\intii}{\int_{-\infty}^\infty}
\newcommand{\la}{\lambda}
\newcommand{\mhf}{^{-1/2}}
\newcommand{\mo}{^{-1}}
\newcommand{\mt}{^{-2}}
\newcommand{\mi}{\mid}
\newcommand{\oom}{{1\over m}}
\newcommand{\oomn}{{1\over mn}}
\newcommand{\oon}{{1\over n}}
\newcommand{\osx}{{{1\over6}}}
\newcommand{\otd}{{{1\over3}}}
\newcommand{\ra}{\to}
\newcommand{\rai}{\ra\infty}
\newcommand{\seight}{_{(8)}}
\newcommand{\sseven}{_{(7)}}
\newcommand{\sfive}{_{(5)}}
\newcommand{\sfour}{_{(4)}}
\newcommand{\ssix}{_{(6)}}
\newcommand{\sk}{_{(k)}}
\newcommand{\snine}{_{(9)}}
\newcommand{\sone}{_{(1)}}

\newcommand{\sthree}{_{(3)}}
\newcommand{\stwo}{_{(2)}}
\newcommand{\sumi}{\sum_i}
\newcommand{\sumiom}{\sum_{i=1}^m}
\newcommand{\sumj}{\sum_j}
\newcommand{\sumjon}{\sum_{j=1}^n}
\newcommand{\thf}{{{1\over2}}}

\newcommand{\var}{\operatorname{var}}
\newcommand{\z}{^0}
\newcommand{\ze}{\zeta}
\newcommand{\seleven}{_{(11)}}
\newcommand{\sten}{_{(10)}}

\newcommand{\lambdahat}{{\widehat\lambda}}
\newcommand{\muhat}{{\widehat\mu}}
\newcommand{\varComp}{\sigma^2}
\newcommand{\varCompmo}{\sigma^{-2}}
\newcommand{\varCompZero}{(\varComp)^0}

\newcommand{\bbeta}{\bolds{\beta}}
\newcommand{\infint}{\int_{-\infty}^{\infty}}
\newcommand{\bbetahat}{\widehat{\bbeta}}
\newcommand{\sigmahat}{{\widehat\sigma}}
\newcommand{\argmaxdum}{\mathop{\arg\max}}
\newcommand{\argmax}[1]{\argmaxdum_{#1}}
\newcommand{\bmu}{\bolds{\mu}}
\newcommand{\smhalf}{{{\frac{1}{2}}}}
\newcommand{\blambda}{\bolds{\lambda}}
\newcommand{\betahat}{{\widehat\beta}}

\newcommand{\uell}{\underline{\ell}}
\newcommand{\muiGVA}{\underline{\muhat}{}_{i}}
\newcommand{\lambdaiGVA}{\underline{\lambdahat}{}_{i}}
\newcommand{\betazGVA}{\underline{\betahat}{}_{0}}
\newcommand{\betaoGVA}{\underline{\betahat}{}_{1}}
\newcommand{\betazZero}{\beta_0^0}
\newcommand{\betaoZero}{\beta_1^0}
\newcommand{\varCompGVA}{{\underline{\sigmahat}}{}^2}

\makeatother

\begin{document}
\begin{frontmatter}

\title{Asymptotic normality and valid inference for
Gaussian variational approximation\thanksref{T1}}
\runtitle{Valid Gaussian variational approximate inference}

\thankstext{T1}{Supported in part by Australian Research
Council grants to the University of Melbourne
and University of Wollongong.}

\begin{aug}
\author[A]{\fnms{Peter} \snm{Hall}\ead[label=e1]{halpstat@ms.unimelb.edu.au}},
\author[B]{\fnms{Tung} \snm{Pham}\ead[label=e2]{tung.pham@epfl.ch}},
\author[C]{\fnms{M. P.} \snm{Wand}\corref{}\ead[label=e3]{matt.wand@uts.edu.au}}
\and
\author[B]{\fnms{S. S. J.} \snm{Wang}\ead[label=e4]{sw918@uow.edu.au}}
\runauthor{Hall, Pham, Wand and Wang}
\affiliation{University of Melbourne, University of Wollongong,
University of Technology, Sydney, and University of Wollongong}
\address[A]{P. Hall\\
Department of Mathematics\\
\quad and Statistics\\
University of Melbourne\\
Melbourne 3000\\
Australia\\
\printead{e1}}
\address[B]{T. Pham\\
S. S. J. Wang\\
Centre for Statistical\\
\quad and Survey Methodology\\
School of Mathematics\\
\quad and Applied Statistics\\
University of Wollongong\\
Wollongong 2522\\
Australia\\
\printead{e2}\\
\hphantom{E-mail: }\printead*{e4}}
\address[C]{M. P. Wand\\
School of Mathematical Sciences\\
University of Technology, Sydney\\
Broadway NSW 2007\\
Australia\\
\printead{e3}} 
\end{aug}

\received{\smonth{11} \syear{2010}}
\revised{\smonth{6} \syear{2011}}

%
\begin{abstract}
We derive the precise asymptotic distributional behavior of Gaussian
variational approximate estimators of the parameters in a
single-predictor Poisson mixed model. These results are the deepest yet
obtained concerning the statistical properties of a variational
approximation method. Moreover, they give rise to asymptotically valid
statistical inference. A simulation study demonstrates that Gaussian
variational approximate confidence intervals possess good to excellent
coverage properties, and have a similar precision to their exact
likelihood counterparts.
\end{abstract}

%
\begin{keyword}[class=AMS]
\kwd[Primary ]{62F12}
\kwd[; secondary ]{62F25}.
\end{keyword}
\begin{keyword}
\kwd{Generalized linear mixed models}
\kwd{longitudinal data analysis}
\kwd{maximum likelihood estimation}
\kwd{Poisson mixed models}.
\end{keyword}

\end{frontmatter}

\section{Introduction}

Variational approximation methods are enjoying an increasing amount of
development and use in statistical problems. This raises questions
regarding their statistical properties, such as consistency of point
estimators and validity of statistical inference. We make significant
inroads into answering such questions via thorough theoretical
treatment of one of the simplest nontrivial settings for which
variational approximation is beneficial: the Poisson mixed model with a
single predictor variable and random intercept. We call this the
\textit{simple Poisson mixed model}.

The model treated here is also treated in \cite{HallOrmerodWand}, but
there attention is confined to bounds and rates of convergence. We
improve upon their results by obtaining the asymptotic distributions of
the estimators. The results reveal that the estimators are
asymptotically normal, have negligible bias and that their variances
decay at least as fast as $m^{-1}$, where $m$ is the number of groups.
For the slope parameter, the faster $(mn)^{-1}$ rate is obtained, where
$n$ is the number of repeated measures.

An important practical ramification of our theory is asymptotically
valid statistical inference for the model parameters. In particular, a
form of studentization leads to theoretically justifiable confidence
intervals for all model parameters. Unlike those based on the exact
likelihood, all Gaussian variational approximate point estimates and
confidence intervals can\vadjust{\goodbreak} be computed without the need for numerical
integration. Simulation results reveal that the confidence intervals
have good to excellent coverage and have about the same length as exact
likelihood-based intervals.

Variational approximation methodology is now a major research area
within computer science; see, for example, Chapter 10 of \cite{Bishop}.
It is beginning to have a presence in statistics as well (e.g.,
\cite{McGroryTitteringtonReevesPettitt,OrmerodWandJCGS}). A summary of
the topic from a statistical perspective is given in
\cite{OrmerodWandTAS}. Late 2008 saw the first beta release of a
software library, \textsf{Infer.NET} \cite{MinkaWinnGuiverKannan}, for
facilitation of variational approximate inference. A high proportion of
variational approximation methodology is framed within Bayesian
hierarchical structures and offers itself as a faster alternative to
Markov chain Monte Carlo methods. The chief driving force is
applications where speed is at a premium and some accuracy can be
sacrificed. Examples of such applications are cluster analysis of
gene-expression data \cite{TeschendorffWangBarbosaMorais}, fitting
spatial models to neuroimage data \cite{FlandinPenny}, image
segmentation \cite{BoccignoneNapoletanoFerraro} and genome-wide
association analysis \cite{LogsdonHoffmanMezey}. Other recent
developments in approximate Bayesian inference include
\textit{approximate Bayesian computing} (e.g.,
\cite{BeaumontZhangBalding}), \textit{expectation propagation} (e.g.,
\cite{Minka}), \textit{integrated nested Laplace approximation} (e.g.,
\cite{RueMartinoChopin}) and \textit{sequential Monte Carlo}
(e.g.,~\cite{DelMoralDoucetJasra}).\looseness=-1

As explained in \cite{Bishop} and \cite{OrmerodWandTAS}, there are many
types of variational approximations. The most popular is
\textit{variational Bayes} (also known as \textit{mean field}
approximation), which relies on product restrictions applied to the
joint posterior densities of a Bayesian model. The present article is
concerned with \textit{Gaussian} \mbox{variational} approximation in
frequentist models containing random effects. There are numerous models
of this general type. One of their hallmarks is the difficulty of exact
likelihood-based inference for the model parameters due to presence of
nonanalytic integrals. Generalized linear mixed models (e.g., Chapter 7
of \cite{McCullochSearleNeuhaus}) form a large class of models for
handling within-group correlation when the response variable is
non-Gaussian. The simple Poisson mixed model lies within this class.
From a theoretical standpoint, the simple Poisson mixed model is
attractive because it possesses the computational challenges that
motivate Gaussian variational approximation---exact likelihood-based
inference requires quadrature---but its simplicity makes it amenable to
deep theoretical treatment. We take advantage of this simplicity to
derive the asymptotic distribution of the Gaussian variational
approximate estimators, although the derivations are still quite
intricate and involved. These results represent the deepest statistical
theory yet obtained for a variational approximation method.

Moreover, for the first time, asymptotically valid inference for a
variational approximation method is manifest. Our theorem reveals that
each estimator is asymptotically normal, centered on the true parameter
value and with a Studentizable variance. Replacement of the unknown
quantities by consistent estimators results in asymptotically valid
confidence intervals and Wald hypothesis tests. A~simulation study
shows that Gaussian variational\vadjust{\goodbreak} approximate confidence intervals
possess good to excellent coverage properties, especially in the case
of the slope parameter.

Section \ref{secspmm} describes the simple Poisson mixed model and
Gaussian variational approximation. An asymptotic normality theorem is
presented in Section \ref{sectheory}. In Section \ref{secinference} we
discuss the implications for valid inference and perform some numerical
evaluations. Section \ref{secproofs} contains the proof of the theorem.

\section{Gaussian variational approximation for the simple
Poisson mixed model}\label{secspmm}

The simple Poisson mixed model that we study here is identical to that
treated in \cite{HallOrmerodWand}. Section 2 of that paper provides a
detailed description of the model and the genesis of Gaussian
variational approximation for estimation of the model parameters. Here
we give just a rudimentary account of the model and estimation
strategy.

The simple Poisson mixed model is
%
%
\begin{eqnarray}\qquad
\label{eqtheModelPtI}
&&Y_{ij}|X_{ij},U_i \mbox{ independent Poisson}\qquad\mbox{with mean }
\exp(\betazZero+\betaoZero X_{ij}+U_i),\\
\label{eqtheModelPtII}
&&U_i \mbox{ independent } N(0,\varCompZero).
\end{eqnarray}
The $X\ourij$ and $U_i$, for $1\le i\le m$ and $1\le j\le n$, are
totally independent random variables, with the $X\ourij$'s distributed
as $X$. We observe values of $(X_{ij},Y_{ij})$, $1\le i\le m$, $1\le
j\le n$, while the $U_i$ are unobserved latent variables. See, for
example, Chapter 7 and Section 14.3 of \cite{McCullochSearleNeuhaus}
for further details on this model and its use in longitudinal data
analysis. In applications it is typically the case that $m\gg n$.

Let $\bbeta\equiv(\beta_0,\beta_1)$ be the vector of fixed effects
parameters.
The conditional log-likelihood of $(\bbeta,\sigma^2)$ is the
logarithm of the
joint probability mass function of the $Y_{ij}$'s, given the $X_{ij}$'s,
as a function of the parameters
%
%
\begin{eqnarray}\label{eqlogLikelihood}
\ell(\bbeta,\sigma^2)&=&\sum_{i=1}^m
\sum_{j=1}^{n}\{Y_{ij}(\beta_0+\beta_1 X_{ij})-\log(Y_{ij}!)\}
- {\frac{m}{2}}\log(2\pi\sigma^2)\nonumber\\[-8pt]\\[-8pt]
&&{}+\sum_{i=1}^m\log\infint\exp\Biggl\{\sum_{j=1}^{n}(Y_{ij}u
-e^{\beta_0+\beta_1 X_{ij}+u})-\frac{u^2}{2\sigma^2}\Biggr\} \,du.
\nonumber
\end{eqnarray}
Maximum likelihood estimation is hindered by the presence
of $m$ intractable integrals in (\ref{eqlogLikelihood}).
However, the $i$th of these integrals can be written as
\begin{eqnarray*}
&&\infint\exp\biggl\{\sum_{j=1}^{n}(Y_{ij}u
-e^{\beta_0+\beta_1 X_{ij}+u})-\frac{u^2}{2\sigma^2}\biggr\}
\frac{e^{-{(1/2)}(u-\mu
_i)^2/\lambda_i}/{\sqrt{2\pi\lambda_i}}}
{e^{-(1/2)(u-\mu_i)^2/\lambda
_i}/{\sqrt{2\pi\lambda_i}}} \,du\\
&&\qquad
=\sqrt{2\pi\lambda_i} E_{ \uU_i}\Biggl[\exp\Biggl\{\sum
_{j=1}^{n}(Y_{ij}\uU_i
-e^{\beta_0+\beta_1 X_{ij}+\uU_i})-\frac{\uU_i^2}{2\sigma^2}
+\frac{(\uU_i-\mu_i)^2}{2\lambda_i}
\Biggr\}\Biggr],
\end{eqnarray*}
where, for $1\le i\le m$, $E_{ \uU_i}$ denotes expectation with respect
to the random variable $\uU_i\sim N(\mu_i,\lambda_i)$
with $\lambda_i>0$.
Jensen's inequality then produces the lower bound
\begin{eqnarray*}
&&\log E_{ \uU_i}\Biggl[\exp\Biggl\{\sum_{j=1}^{n}(Y_{ij}\uU_i
-e^{\beta_0+\beta_1 X_{ij}+\uU_i})-\frac{\uU_i^2}{2\sigma^2}
+\frac{(\uU_i-\mu_i)^2}{2\lambda_i}
\Biggr\}\Biggr]\\
&&\qquad\ge E_{ \uU_i}\Biggl\{\sum_{j=1}^{n}(Y_{ij}\uU_i
-e^{\beta_0+\beta_1 X_{ij}+\uU_i})-\frac{\uU_i^2}{2\sigma^2}
+\frac{(\uU_i-\mu_i)^2}{2\lambda_i}
\Biggr\},
\end{eqnarray*}
which is tractable. Standard manipulations then lead to
%
%
\begin{equation}\label{eqlogLikLB}
\ell(\bbeta,\sigma^2) \ge\underline{\ell}(\bbeta,\sigma
^2,\bmu,\blambda)
\end{equation}
for all vectors $\bmu=(\mu_1,\ldots,\mu_m)$
and $\blambda=(\lambda_1,\ldots,\lambda_m)$, where
%
%
\begin{eqnarray}\label{eqdiagLambda}
&&\underline{\ell}(\bbeta,\sigma^2,\bmu,\blambda)\nonumber\\
&&\qquad\equiv
\sum_{i=1}^m\sum_{j=1}^n\{Y_{ij}(\beta_0+\beta_1 X_{ij}+\mu_i)
-e^{\beta_0+\beta_1 X_{ij}+\mu_i+\lambda_i/2}
-\log(Y_{ij}!)\}\nonumber\\[-8pt]\\[-8pt]
&&\qquad\quad{} -\frac{m}{2}\log(\sigma^2)+\frac{m}{2}
-\frac{1}{2\sigma^2}\sum_{i=1}^m(\mu_i^2+\lambda_i)\nonumber\\
&&\qquad\quad{}
+\smhalf\sum_{i=1}^m\log(\lambda_i)
\nonumber
\end{eqnarray}
is a \textit{Gaussian variational approximation} to $\ell(\bbeta
,\sigma^2)$.
The vectors $\bmu$ and $\blambda$ are \textit{variational parameters}
and should be chosen to make
$\underline{\ell}(\bbeta,\sigma^2,\bmu,\blambda)$
as close as possible to $\ell(\bbeta,\sigma^2)$.
In view of (\ref{eqlogLikLB}) the Gaussian variational
approximate maximum likelihood estimators are naturally
defined to be
\[
(\underline{\bbetahat},\underline{\sigmahat}^2)
=\mbox{$(\bbeta,\sigma^2)$ component of }
\argmax{\bbeta,\sigma^2,\bmu,\blambda}\uell(\bbeta,\sigma
^2,\bmu,\blambda).
\]

\section{Asymptotic normality results}\label{sectheory}


Consider random variables $(X_{ij},Y_{ij},U_i)$
satisfying (\ref{eqtheModelPtI}) and (\ref{eqtheModelPtII}).
Put
\[
Y\ib=\sum_{i=1}^n Y\ourij\quad\mbox{and}\quad
B_i=\sum_{j=1}^n\exp(\be_0+\be_1 X\ourij),
\]
and consider the following decompositions of the
exact log-likelihood and its Gaussian variational approximation:
\begin{eqnarray*}
\ell(\bbeta,\varComp)&=&\ell_0(\bbeta,\varComp)
+\ell_1(\bbeta,\varComp)+\mbox{DATA} ,\\
\uell(\bbeta,\varComp,\bmu,\blambda)&=&\ell_0(\bbeta,\varComp)
+\ell_2(\bbeta,\varComp,\bmu,\blambda)+\mbox{DATA} ,
\end{eqnarray*}
where
%
%
\begin{eqnarray}
\label{eqonePtOne}
\ell_0(\bbeta,\varComp)&=&\sumiom\sumjon Y\ourij(\be_0+\be_1
X\ourij)
-\thf m \log\varComp,\nonumber\\[-6pt]\\[-6pt]
\ell_1(\bbeta,\varComp)&=&\sumiom\log
\biggl\{\intii\exp\biggl(Y\ib u-B_i e^u-\thf\varCompmo u^2\biggr)
\,du\biggr\} ,\nonumber\\[3pt]
\label{eqonePtTwo}
\ell_2(\bbeta,\varComp,\bmu,\blambda)&=&\sumiom\biggl\{\mu_i
Y\ib
-B_i \exp\biggl(\mu_i+\thf\la_i\biggr)\biggr\}\nonumber\\[-6pt]\\[-6pt]
&&{} -\thf\varCompmo\sumiom(\mu
_i^2+\la_i)
+\thf\sumiom\log\la_i ,\nonumber
\end{eqnarray}
and DATA denotes a quantity depending on the $Y_{ij}$ alone,
and not on $\bbeta$ or $\varComp$. Note that
\[
\uell(\bbeta,\varComp)=\max_{\bmu,\blambda} \uell(\bbeta
,\varComp,\bmu,\blambda)
=\uell_0(\bbeta,\varComp)+\max_{\bmu,\blambda}
\uell_2(\bbeta,\varComp,\bmu,\blambda).
\]

Our upcoming theorem relies on the following assumptions:
\begin{longlist}[(A2)]
\item[(A1)] the moment generating function of $X$, $\phi(t)=E\{\exp
(tX)\}$,
is well defined on the whole real line;
\item[(A2)] the mapping that takes $\be$ to $\phi'(\be)/\phi(\be
)$ is invertible;
\item[(A3)] in some neighborhood of $\be_1\z$ (the true value of
$\be_1$),
$(d^2/d\be^2) \log\phi(\be)$ does not vanish;
\item[(A4)] $m=m(n)$ diverges to infinity with $n$, such that $n/m\ra0$
as $n\rai$;
\item[(A5)] for a constant $C>0$, $m=O(n^C)$ as $m$ and $n$ diverge.
\end{longlist}

Define
%
%
\begin{equation}\label{eqonePtFive}
\tau^2=\frac{\exp\{-\varCompZero/2-\be_0\z\} \phi
(\be_1\z)}
{\phi''(\be_1\z)\phi(\be_1\z)-\phi'(\be_1\z)^2} .
\end{equation}
%
\noindent
The precise asymptotic behavior of $\betazGVA$, $\betaoGVA$ and
$\varCompGVA$
is conveyed by:
\begin{theorem}
\label{thmmainThm}
Assume that conditions \textup{(A1)--(A5)} hold. Then
%
%
\begin{equation}\label{eqPGHa}
\betazGVA-\betazZero=m\mhf N_0+o_p(n\mo+m\mhf),
\end{equation}
where the random variable $N_0$ is normal $N(0,\varCompZero)$;
%
%
\begin{equation}\label{eqPGHb}
\betaoGVA-\betaoZero=(mn)\mhf N_1+o_p\{n\mt+(mn)\mhf\},
\end{equation}
where the random variable $N_1$
is normal $N(0,\tau^2)$; and
%
%
\begin{equation}\label{eqPGHc}
\varCompGVA-\varCompZero=m\mhf N_2+o_p(n\mo+m\mhf),
\end{equation}
where the random variable $N_2$
is normal $N(0,2\{\varCompZero\}^2)$.
\end{theorem}
\begin{remark*}
All three Gaussian variational approximate estimators have
asymptotically normal distributions with asymptotically
negligible bias. The estimators $\betazGVA$
and $\varCompGVA$ have variances of size $m^{-1}$, as $m$
and $n$ diverge in such a manner that $n/m\ra0$.
The estimator $\betaoGVA$ has variance of size $(mn)^{-1}$.
Hence, the estimator $\betaoGVA$ is distinctly more accurate than
either $\betazGVA$ or $\varCompGVA$, since it converges to the respective
true parameter value at a strictly faster rate.
For the estimator $\betaoGVA$, increasing both
$m$ and $n$ reduces variance. However, in the cases of the
estimators $\betazGVA$ or $\varCompGVA$,
only an increase in $m$ reduces variance.
\end{remark*}

\section{Asymptotically valid inference}\label{secinference}

Theorem \ref{thmmainThm} reveals that $\betazGVA$, $\betaoGVA$ and
$\varCompGVA$
are each asymptotically normal with means corresponding to the true parameter
values. The variances depend on known functions of the parameters and
$\phi(\betaoZero)$,
$\phi'(\betaoZero)$ and $\phi''(\betaoZero)$. Since the latter
three quantities
can be estimated unbiasedly via
\begin{eqnarray*}
\widehat{\phi(\betaoZero)}&=&\frac{1}{mn}\sum_{i=1}^m\sum
_{j=1}^n\exp(X_{ij}\betaoGVA),\\
\widehat{\phi'(\betaoZero)}&=&\frac{1}{mn}\sum_{i=1}^m\sum
_{j=1}^n X_{ij} \exp(X_{ij}\betaoGVA)
\end{eqnarray*}
and
\[
\widehat{\phi''(\betaoZero)}=\frac{1}{mn}\sum_{i=1}^m\sum
_{j=1}^n X_{ij}^2 \exp(X_{ij}\betaoGVA),
\]
we can consistently estimate the asymptotic variances for inferential procedures
such as confidence intervals and Wald hypothesis tests. For example,
the quantity $\tau^2$ appearing in the expression for the
asymptotic variance of $\betaoGVA$ can be consistently estimated
by\vadjust{\eject}
\[
\widehat{\tau}^2=\frac{\exp(-\varCompGVA/2-\betazGVA
)
\widehat{\phi(\be_1\z)}}
{\widehat{\phi''(\be_1\z)}\widehat{\phi(\be_1\z)}-\widehat
{\phi'(\be_1\z)}{}^2}.
\]
Approximate $100(1-\alpha)\%$ confidence intervals for
$\betazZero$, $\betaoZero$
and $\varCompZero$ are
%
%
\begin{eqnarray}\label{eqCIs}
&\displaystyle \betazGVA\pm\Phi\biggl(1-\smhalf\alpha\biggr)\sqrt{\frac{\varCompGVA
}{m}},\qquad
\betaoGVA\pm\Phi\biggl(1-\smhalf\alpha\biggr)\sqrt{\frac{\widehat{\tau}^2}{mn}}
\quad\mbox{and}&\nonumber\\[-8pt]\\[-8pt]
&\displaystyle \varCompGVA\pm\Phi\biggl(1-\smhalf\alpha\biggr)
\varCompGVA\sqrt{\frac{2}{m}},&\nonumber
\end{eqnarray}
where $\Phi$ denotes the $N(0,1)$ distribution function.
These confidence intervals are \textit{asymptotically valid} since
they involve studentization based on consistent estimators
of all unknown quantities.

We ran a simulation study to evaluate the coverage properties
of the Gaussian variational approximate confidence intervals
(\ref{eqCIs}). The true parameter vector
$(\betazZero,\betaoZero,\varCompZero)$ was
allowed to vary over
\begin{eqnarray*}
&&\{(-0.3,0.2,0.5),(2.2,-0.1,0.16),\\
&&\qquad(1.2,0.4,0.1),(0.02,1.3,1),
(-0.3,0.2,0.1)\},
\end{eqnarray*}
and the distribution of the $X_{ij}$ was taken to be either
$N(0,1)$ or $\mbox{Uniform}(-1,1)$, the uniform distribution
over the interval $(-1,1)$. The number groups $m$ varied
over $100,200,\ldots,1\mbox{,}000$ with $n$ fixed at $m/10$
throughout the study. For each of the ten possible
combinations of true parameter vector and $X_{ij}$ distribution,
and sample size pairs, we generated 1,000 samples and
computed 95\% confidence intervals based on (\ref{eqCIs}).

\begin{sidewaysfigure}

\includegraphics{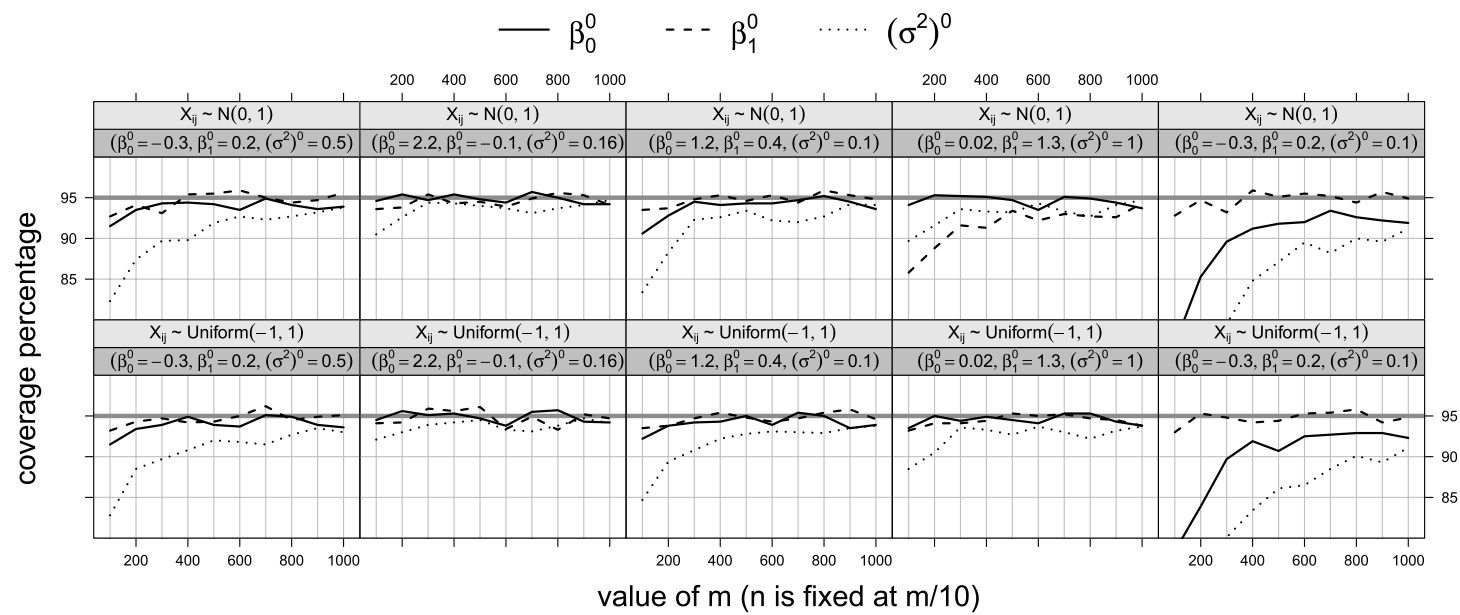}

\caption{Actual coverage percentage of nominally 95\% Gaussian
variational approximate confidence intervals for the
parameters in the simple Poisson mixed model.
The nominal percentage is shown as a thick grey horizontal line.
The percentages are based on 1,000 replications.
The values of $m$ are
$100,200,\ldots,1\mbox{,}000$. The value of $n$ is fixed at $n=m/10$.}
\label{figCIsimPlot}
\end{sidewaysfigure}

Figure \ref{figCIsimPlot} shows the actual coverage
percentages for the nominally 95\% confidence intervals.
In the case of $\betaoZero$, the actual and nominal
percentages are seen to have very good agreement, even for
$(m,n)=(100,10)$. This is also the case for $\betazZero$
for the first four true parameter vectors. For the fifth
one, which has a relatively low amount of within-subject
correlation, the asymptotics take a bit longer to become
apparent, and we see that $m\ge400$ is required to get the
actual coverage above 90\%, that is, within 5\% of the nominal
level. For $\varCompZero$, a similar comment applies,
but with $m\ge800$. The superior coverage
of the $\betaoZero$ confidence intervals is in keeping
with the faster convergence rate apparent from Theorem \ref{thmmainThm}.

%

Lastly, we ran a smaller simulation study to check
whether or not the lengths of the Gaussian variational approximate
confidence intervals are compromised in achieving the
good coverage apparent in Figure \ref{figCIsimPlot}.
For each of the same settings used to produce that figure
we generated 100 samples and computed the exact
likelihood-based confidence intervals using
adaptive Gauss--Hermite quadrature (via the \textsf{R} language
\cite{Rteam} package \texttt{lme4} \cite{BatesMaechler}).
In almost every case, the Gaussian variational approximate
confidence intervals were slightly shorter than their
exact counterparts. This reassuring result indicates that
the good coverage performance is not accompanied by a
decrease in precision.

\section{\texorpdfstring{Proof of Theorem \protect\ref{thmmainThm}}{Proof of Theorem 3.1}}
\label{secproofs}

The proof Theorem \ref{thmmainThm} requires some additional
notation, as well as several stages of asymptotic approximation.
This section provides full details, beginning with definitions
of the necessary notation.

\subsection{Notation}

Recall that $\be_0\z$, $\be_1\z$ and $\varCompZero$ denote the
true values of parameters and that $\betazGVA$, $\betaoGVA$ and
$\varCompGVA$ denote their respective Gaussian variational
approximate estimators.

%
\begin{table}
\tablewidth=310pt
\caption{Definitions of the $O_{(k)}$ notation used in the proofs}
\label{tabOnotation}
\begin{tabular*}{\tablewidth}{@{\extracolsep{\fill}}lc@{}}
\hline
\textbf{Notation} & \textbf{Meaning} \\
\hline
$O\sone$ & $O_p(m\mhf+n\mo)$ \\
$O\stwo$ &$O_p(m\mo+n\mt)$\\
$O\sthree$&$O(n^{\ep-(1/2)})$,
uniformly in $1\le i\le m$, for each $\ep>0$\\
$O\sfour$ &$O(n^{\ep-1})$,
uniformly in $1\le i\le m$, for each $\ep>0$\\
$O\sfive$ &$O(n^{\ep-(3/2)})$,
uniformly in $1\le i\le m$, for each $\ep>0$\\
$O\ssix$&$O_p(m\mo+n^{\ep-(3/2)})$,
uniformly in $1\le i\le m$, for each $\ep>0$\\
$O\sseven$ &$O_p\{(m\mo+n\mt) n^{\ep-(1/2)}\}$,
uniformly in $1\le i\le m$, for each $\ep>0$\\
$O\seight$ &$O_p\{(m\mhf+n\mo)^3 n^\ep\}$,
uniformly in $1\le i\le m$, for each $\ep>0$\\
$O\snine$&$O_p\{(mn)^{-1/2}+n^{\ep-(3/2)}\}$,
uniformly in $1\le i\le m$, for each $\ep>0$\\
$O\sten$&$O_p\{(m\mhf+n^{-5/2}) n^\ep\}$,
uniformly in $1\le i\le m$, for each $\ep>0$\\
$O\seleven$&$O_p\{(m\mhf n\mo+n\mt) n^\ep\}$,
uniformly in $1\le i\le m$, for each $\ep>0$\\
\hline
\end{tabular*}
\end{table}

The proofs use ``$O\sk$'' notation, for $k=1,\ldots,11$,
as defined in Table \ref{tabOnotation}.

\subsection{Formulae for estimators}

First we give, in (\ref{eqTwoPtTwo})--(\ref{eqTwoPtSix}) below, the results
of equating to zero the
derivatives of $\ell_0(\be,\varComp)+\ell_2(\be,\varComp,\la
,\mu)$ with respect
to $\be_0$, $\be_1$, $\varComp$, $\lambda_i$ and $\mu_i$, respectively:
%
%
\begin{eqnarray}
\label{eqTwoPtTwo}
\sumiom\biggl\{Y\ib-B_i \exp\biggl(\muiGVA+\thf\lambdaiGVA
\biggr)\biggr\}
&=&0 ,\\
\label{eqTwoPtThree}
\qquad\sumiom\sumjon X\ourij\biggl\{Y\ourij
-\exp\biggl(\betazGVA+\muiGVA+\thf\lambdaiGVA+\betaoGVA X\ourij
\biggr) \biggr\}
&=&0 ,\\
\label{eqTwoPtFour}
\oom\sumiom(\lambdaiGVA+\muiGVA^2)
&=&\varCompGVA
,\\
\label{eqTwoPtFive}
\lambdaiGVA\mo-B_i \exp\bigl(\muiGVA+\tfrac{1}{2}\lambdaiGVA
\bigr)-(\varCompGVA)^{-1}
&=&0 ,\qquad 1\le i\le m,\\
\label{eqTwoPtSix}
Y\ib-B_i \exp\bigl(\muiGVA+\tfrac{1}{2}\lambdaiGVA\bigr)-(\varCompGVA)^{-1}
\muiGVA&=&0 ,\qquad 1\le i\le m.
\end{eqnarray}
These are the analogs of the likelihood equations in the
conventional approach to inference.

The next step is to put (\ref{eqTwoPtTwo}), (\ref{eqTwoPtThree})
and (\ref{eqTwoPtSix})
into more accessible form, in (\ref{eqTwoPtSeven}),
(\ref{eqTwoPtTwelve}) and (\ref{eqTwoPtThirteen}), respectively.
Adding (\ref{eqTwoPtSix}) over $1\le i\le m$ and subtracting the result
from (\ref{eqTwoPtTwo}) we deduce that
%
%
\begin{equation}\label{eqTwoPtSeven}
\sumiom\muiGVA=0 .
\end{equation}
Defining
\[
\De=\oomn\sumiom\sumjon X\ourij
\{Y\ourij-\exp(\be_0\z+\be_1\z X\ourij+U_i)
\}
\]
we deduce that (\ref{eqTwoPtThree}) is equivalent to
%
%
\begin{eqnarray}\label{eqTwoPtEight}
&&\De+\exp(\be_0\z)
\oomn\sumiom\sumjon X\ourij\exp(U_i+\be_1\z
X\ourij)\nonumber\\[-8pt]\\[-8pt]
&&\qquad{} -\exp(\be_0) \oomn\sumiom\sumjon X\ourij
\exp\biggl(\muiGVA+\thf\lambdaiGVA+\be_1 X\ourij\biggr)=0 .
\nonumber
\end{eqnarray}

Define $\xi_i$, $\eta_i$ and $\ze_i$ by, respectively,
%
%
\begin{eqnarray}
\label{eqTwoPtNine}
\oon\sumjon X\ourij\exp(\be_1\z X\ourij)
&=&\phi'(\be_1\z)
\exp(\xi_i) ,\\
\label{eqTwoPtTen}
\oon\sumjon X\ourij
\exp(\betaoGVA X\ourij)
&=&\phi'(\betaoGVA)
\exp(\eta_i) ,\\[-14pt]\nonumber
\end{eqnarray}
%
%
\begin{eqnarray}\label{eqTwoPtEleven}
&&\exp\biggl(\betazGVA+\muiGVA+\thf\lambdaiGVA\biggr)
\oon\sumjon\{\exp(\betaoGVA X\ourij)
-\phi(\betaoGVA)\}
\nonumber\\
&&\qquad=
\exp(\be_0\z+U_i)\Biggl[\phi(\be_1\z) \{
1-\exp(\ze_i)\}
\nonumber\\[-8pt]\\[-8pt]
&&\qquad\quad\hphantom{\exp(\be_0\z+U_i)\Biggl[}
{}
+\oon\sumjon\{Y\ourij\exp(-\be_0\z-U_i)
-\phi(\be_1\z)\}\Biggr]\nonumber\\
&&\qquad\quad{}
-(\varCompGVA n)\mo\muiGVA.
\nonumber
\end{eqnarray}

With probability converging to 1 as $n\rai$ the definitions at
(\ref{eqTwoPtNine})--(\ref{eqTwoPtEleven}) are
valid simultaneously for all $1\le i\le m$, because the variables $\xi
_i$, $\eta_i$
and $\ze_i$ so defined converge to zero, uniformly in $1\le i\le m$,
in probability.
See (\ref{eqTwoPtThirtyOne}), (\ref{eqTwoPtThirtyTwo}) and (\ref
{eqTwoPtTwentySix})
below for approximations
to $\xi_i$, $\eta_i$ and $\ze_i$;
indeed, those formulae quickly imply that each of $\xi_i$, $\eta_i$
and $\ze_i$ equals
$O\sthree$.

Without loss of generality, $\phi'(t)$ is bounded away from zero in a
neighborhood
of $\be_1\z$. Indeed, if the latter property does not hold, simply
add a constant
to the random variable $X$ to ensure that $\phi'(\be_1\z)\neq0$. We
assume that
$\be_1\z$ is in the just-mentioned neighborhood, and we consider only
realizations
for which $\be_1$ is also in the neighborhood. (The latter property
holds true
with probability converging to 1 as $n\rai$.) The definition of $\ze
_i$ at
(\ref{eqTwoPtEleven})
can be justified using the fact that $\muiGVA<Y\ib$, as shown in
Theorem 2 of
\cite{HallOrmerodWand}.

In this notation we can write (\ref{eqTwoPtEight}) as
%
%
\begin{eqnarray}\label{eqTwoPtTwelve}
&&\De+\phi'(\be_1\z)
\oom\sumiom\exp(\be_0\z+U_i+\xi_i)
\nonumber\\[-8pt]\\[-8pt]
&&\qquad
=\phi'(\betaoGVA)
\oom\sumiom\exp\biggl(\betazGVA+\muiGVA+\thf\lambdaiGVA+\eta
_i\biggr)
\nonumber
\end{eqnarray}
and write (\ref{eqTwoPtSix}) as
%
%
\begin{equation}\label{eqTwoPtThirteen}
\exp\bigl(\betazGVA+\muiGVA+\tfrac{1}{2}\lambdaiGVA\bigr) \phi(\betaoGVA)
=\exp(\be_0\z+U_i+\ze_i) \phi(\be_1\z) .
\end{equation}

Substituting (\ref{eqTwoPtThirteen}) into (\ref{eqTwoPtTwelve}) we obtain
%
%
\begin{eqnarray}\label{eqTwoPtFourteen}
&&
\De\exp(-\be_0\z) \phi(\be_1\z)\mo
+\phi'(\be_1\z) \phi(\be_1\z)\mo
\oom\sumiom\exp(U_i+\xi_i)\nonumber\\[-8pt]\\[-8pt]
&&\qquad=\phi'(\betaoGVA) \phi(\betaoGVA)\mo
\oom\sumiom\exp(U_i+\eta_i+\ze_i) .
\nonumber
\end{eqnarray}

\subsection{\texorpdfstring{Approximate formulae for $U_i$ and
$\lambdaiGVA$}{Approximate formulae for $U_i$ and $lambdai_i$}}

The formulae are given at (\ref{eqTwoPtSeventeen}) and (\ref
{eqTwoPtNineteen}), respectively. To derive them,
note that (\ref{eqTwoPtSix}) implies that
\begin{eqnarray*}
&&\bigl(1+O\sthree\bigr) \phi(\be_1\z) \exp(\be_0\z
+U_i)\\
&&\qquad{} -\bigl(1+O\sthree\bigr) \phi(\be_1\z)
\exp\bigl(\be_0\z+\muiGVA+\tfrac{1}{2}\lambdaiGVA\bigr)
-(n\varCompGVA)\mo\muiGVA=0 .
\end{eqnarray*}
Here we have used the fact that, by \cite{HallOrmerodWand},
%
%
\begin{equation}\label{eqTwoPtFifteen}
\betazGVA-\be_0\z=O\sone,\qquad \betaoGVA-\be_1\z=O\sone,
\end{equation}
and that by (1.3), $\max_{1\le i\le m} |X_i|=O_p(n^\ep)$ for all
$\ep>0$.
Therefore,
%
%
\begin{equation}\label{eqTwoPtSixteen}\quad
\bigl(1+O\sthree\bigr) \exp(U_i)
=\bigl(1+O\sthree\bigr) \exp\bigl(\muiGVA+\tfrac{1}{2}\lambdaiGVA\bigr)
+(cn\varCompGVA)\mo\muiGVA,
\end{equation}
where $c=\phi(\be_1\z) \exp(\be_0\z)$. The result
$\max_{1\le i\le m} |U_i|=O_p\{(\log n)\half\}$ follows from
properties of
extrema of Gaussian variables and the fact that $m=O(n^C)$ for a
constant $C>0$.
Moreover, by Theorem 2 of \cite{HallOrmerodWand}, $0<\lambdaiGVA
<\varCompGVA$.
Therefore (\ref{eqTwoPtSixteen}) implies that
$\max_{1\le i\le n} |\muiGVA|=O_p\{(\log n)\half\}$. [Note that,
for any
constant $C>0$, $\exp\{-C (\log n)\half\}=n^{-C (\log n)\mhf}$,
which is of
larger order than $n^{-\ep}$ for each $\ep>0$.]
Hence, by (\ref{eqTwoPtSixteen}),
\[
\bigl(1+O\sthree\bigr) \exp(U_i) =\bigl(1+O\sthree\bigr)
\exp\bigl(\muiGVA+\tfrac{1}{2}\lambdaiGVA\bigr) ,
\]
and so, taking logarithms,
%
%
\begin{equation}\label{eqTwoPtSeventeen}
U_i=\muiGVA+\tfrac{1}{2}\lambdaiGVA+O\sthree.
\end{equation}

Formula (\ref{eqTwoPtFive}) and property (\ref{eqTwoPtFifteen}) entail
%
%
\begin{equation}\label{eqTwoPtEighteen}
(n\lambdaiGVA)\mo
-\bigl(1+O\sthree\bigr) \phi(\be_1\z)
\exp\bigl(\muiGVA+\tfrac{1}{2}\lambdaiGVA+\be_0\z\bigr)
-(n\varCompGVA)\mo=0 .
\end{equation}
Using (\ref{eqTwoPtSeventeen}) to substitute
$U_i+O\sthree$ for $\muiGVA+\thf\lambdaiGVA$
in (\ref{eqTwoPtEighteen}) we
deduce from that result that
\begin{eqnarray*}
(n\lambdaiGVA)\mo&=&\bigl(1+O\sthree\bigr) \phi(\be_1\z)
\exp(U_i+\be_0\z)
+(n\varCompGVA)\mo\\
&=&\bigl(1+O\sthree\bigr) \phi(\be_1\z)
\exp(U_i+\be_0\z) ,
\end{eqnarray*}
where to obtain the second identity we again used the fact that
\[
\max_{1\le i\le m} |U_i|=O_p\{(\log n)\half\}.
\]
Therefore,
%
%
\begin{eqnarray}\label{eqTwoPtNineteen}
\lambdaiGVA&=&\bigl(1+O\sthree\bigr) \{n \phi(\be_1\z)
\exp(U_i+\be_0\z)\}\mo\nonumber\\[-8pt]\\[-8pt]
&=&\{n \phi(\be_1\z)
\exp(U_i+\be_0\z)\}\mo+O\sfive,
\nonumber
\end{eqnarray}
where $O\sfive$ is as defined in Table \ref{tabOnotation}.
To obtain the second identity in
(\ref{eqTwoPtNineteen}) we used the fact that
$\max_{1\le i\le m} \exp(-U_i)=O(n^\ep)$ for all $\ep>0$.

\subsection{\texorpdfstring{Initial approximations to $\betazGVA-\be_0\z$ and
$\betaoGVA-\be_1\z$}{Initial approximations to $beta_0- beta_0^0$ and
$beta _1-beta_1^0$}}

These approximations are given at (\ref{eqTwoPtTwenty}), (\ref
{eqTwoPtTwentyTwo})
and (\ref{eqTwoPtThirty}), and lead to central
limit theorems for $\betaoGVA-\be_1\z$, $\betazGVA-\be_0\z$ and
$\varCompGVA-\varCompZero$, respectively.
To derive the approximations, write $\ga(\be_1)=\phi'(\be_1) \phi
(\be_1)\mo$ and
note that, defining $O\stwo$ as in Table \ref{tabOnotation}, we have
\begin{eqnarray*}
\ga(\betaoGVA)&=&\ga(\be_1\z)+(\betaoGVA-\be_1\z)
\ga'(\be_1\z)
+O_p(|\betaoGVA-\be_1\z|^2)\\
&=&\ga(\be_1\z)+\{1+O_p(m^{-1/2}+n\mo)
\}
(\betaoGVA-\be_1\z) \ga'(\be_1\z) .
\end{eqnarray*}
[Here we have used (\ref{eqTwoPtFifteen}).] Therefore, by (\ref
{eqTwoPtFourteen})
and for each $\ep>0$,
\begin{eqnarray*}
&&\De\exp(-\be_0\z) \phi(\be_1\z)\mo
+\ga(\be_1\z)
\oom\sumiom\exp(U_i+\xi_i)\\
&&\qquad=[\ga(\be_1\z)
+\{1+O_p(m^{-1/2}+n\mo)\}
(\betaoGVA-\be_1\z) \ga'(\be_1\z)]\\
&&\qquad\quad{} \times
\oom\sumiom\exp(U_i+\eta_i+\ze_i) .
\end{eqnarray*}
That is,
%
%
\begin{eqnarray}\label{eqTwoPtTwenty}
&&(\betaoGVA-\be_1\z) \ga'(\be_1\z)
\oom\sumiom\exp(U_i+\eta_i+\ze_i)\nonumber\\
&&\qquad=\ga(\be_1\z)
\oom\sumiom\exp(U_i) \{\exp(\xi_i)
-\exp(\eta_i+\ze_i)\}\\
&&\qquad\quad{} +\De\exp(-\be_0\z) \phi(\be_1\z
)\mo
+O\stwo.
\nonumber
\end{eqnarray}

Taking logarithms of both sides of (\ref{eqTwoPtThirteen}) we obtain
%
%
\begin{equation}\label{eqTwoPtTwentyOne}
\log\{\phi(\betaoGVA)/\phi(\be_1\z)\}
=\be_0\z-\betazGVA+U_i+\ze_i-\muiGVA-\tfrac{1}{2}\lambdaiGVA,
\end{equation}
which, on adding over $i$ and dividing by $m$, implies that
\[
\log\{\phi(\betaoGVA)/\phi(\be_1\z)\}
=\be_0\z-\betazGVA+\oom\sumiom\biggl(U_i+\ze_i-\muiGVA-\thf
\lambdaiGVA\biggr) ,
\]
which in turn gives
%
%
\begin{eqnarray}\label{eqTwoPtTwentyTwo}\quad
\betazGVA-\be_0\z&=&-(\betaoGVA-\be_1\z) \ga(\be
_1\z)
+\oom\sumiom\biggl(U_i+\ze_i-\muiGVA-\thf\lambdaiGVA\biggr)
+O\stwo\nonumber\\
&=&-(\betaoGVA-\be_1\z) \ga(\be_1\z)
+\oom\sumiom(U_i+\ze_i)\\
&&{}
-\biggl\{2n \phi(\be_1\z)
\exp\biggl(\be_0\z-\thf\varCompZero\biggr)\biggr\}\mo
+O\ssix,
\nonumber
\end{eqnarray}
where we used (\ref{eqTwoPtNineteen}) to substitute for $\lambdaiGVA
$ and (\ref{eqTwoPtSeven}) to eliminate $\muiGVA$
from the right-hand side, and employed (\ref{eqTwoPtFifteen}) to
bound $(\betaoGVA-\be_1\z)^2$.
Note too that $E\{\exp(-U_i)\}=\exp(\thf\varCompZero)$; a term
involving $E\{\exp(-U_i)\}$
arises from $\sumi\lambdaiGVA$ via (\ref{eqTwoPtNineteen}).

\subsection{\texorpdfstring{Approximation to $\ze_i$}{Approximation to $zeta_i$}}

The approximation is given at (\ref{eqTwoPtTwentySix}). First we derive
an expansion, at (\ref{eqTwoPtTwentyThree}) below, of $\muiGVA$.
Reflecting (\ref{eqTwoPtSeventeen}), define the random variable
$\delta_i$ by $\muiGVA=U_i-\thf\lambdaiGVA+\delta_i$.
Then,\vspace*{1pt} by (\ref {eqTwoPtSeventeen}), $\delta_i=O\sthree$.
Define too $B_{ik}\z=\sumj X\ourij^k \exp(\be_0\z+\be_1\z X\ourij)$ for
$k=0,1,2$, and $\De_i=Y\ib-B_{i0}\z\exp(U_i)$;\vspace*{1pt} and let
$\cF_i$ denote the sigma-field generated by $U_i$ and
$X_{i1},\ldots,X_{in}$. Then $E(\De_i\mi\cF _i)=0$ and
\begin{eqnarray*}
B_i&=&\bigl\{1+\betazGVA-\be_0\z+\tfrac{1}{2}(\betazGVA-\be_0\z)^2
\bigr\} B_{i0}\z\\[3pt]
&&{}+\{\betaoGVA-\be_1\z
 +(\betazGVA-\be_0\z) (\betaoGVA-\be_1\z
)\} B_{i1}\z\\[3pt]
&&{} +\tfrac{1}{2}(\betaoGVA-\be_1\z)^2 B_{i2}
+O\seight,
\end{eqnarray*}
uniformly in $1\le i\le m$ for each $\ep>0$, where $O\seight$ is as in
Table \ref{tabOnotation}.
Therefore,
\begin{eqnarray*}
&&Y\ib-B_i \exp(U_i+\delta_i)\\[3pt]
&&\qquad=Y\ib
-\bigl[\bigl\{1+\betazGVA-\be_0\z+\tfrac{1}{2}(\betazGVA-\be_0\z
)^2\bigr\} B_{i0}\z\\[3pt]
&&\qquad\quad\hspace*{30pt}{}
+\{\betaoGVA-\be_1\z+(\betazGVA-\be_0\z)
(\betaoGVA-\be_1\z)\} B_{i1}\z
\\[3pt]
&&\qquad\quad\hspace*{120pt}{}
+\tfrac{1}{2}(\betaoGVA-\be_1\z)^2 B_{i2}\z\bigr]\\[3pt]
&&\hspace*{58pt}{}\times\exp(U_i)
\bigl(1+\delta_i+\tfrac{1}{2}\delta_i^2+O\sfive\bigr)
+n O\seight,
\end{eqnarray*}
where $O\sfive$ is as in Table \ref{tabOnotation}.
Therefore, defining
\begin{eqnarray*}
\chi_i&=&\bigl\{\betazGVA-\be_0\z+\tfrac{1}{2}(\betazGVA-\be_0\z
)^2\bigr\} B_{i0}\z
+\{\betaoGVA-\be_1\z+(\betazGVA-\be_0\z)
(\betaoGVA-\be_1\z)\} B_{i1}\z\\[3pt]
&&{} +\tfrac{1}{2}(\betaoGVA-\be_1\z)^2 B_{i2}\z,
\end{eqnarray*}
we see that the left-hand side of (\ref{eqTwoPtSix}) equals
\begin{eqnarray*}
&&Y\ib-B_i \exp(U_i+\delta_i)-(\varCompGVA)^{-1} \muiGVA\\
&&\qquad=\De_i-B_{i0}\z\exp(U_i) \bigl(\delta_i+\tfrac{1}{2}
\delta_i^2+O\sfive\bigr)\\
&&\qquad\quad{}
-\chi_i \exp(U_i) \bigl(1+\delta_i+\tfrac{1}{2}\delta_i^2+O\sfive
\bigr)\\
&&\qquad\quad{}
-(\varCompGVA)^{-1} \bigl(U_i-\tfrac{1}{2}\lambdaiGVA+\delta_i\bigr)+n
O\seight\\
&&\qquad=\De_i-\bigl\{\chi_i \exp(U_i)+(\varCompGVA)^{-1}
\bigl(U_i-\tfrac{1}{2}\lambdaiGVA\bigr)\bigr\}\\
&&\qquad\quad{}
-\delta_i \{(B_{i0}\z+\chi_i) \exp
(U_i)+(\varCompGVA)^{-1}\}\\
&&\qquad\quad{} -\tfrac{1}{2}\delta_i^2 (B_{i0}\z
+\chi_i) \exp(U_i)
+n O\sfive+n O\seight.
\end{eqnarray*}
Hence, (\ref{eqTwoPtSix}) implies that
\begin{eqnarray*}
&&\delta_i+\thf\delta_i^2 {(B_{i0}\z+\chi_i) \exp(U_i)
\over(B_{i0}\z+\chi_i) \exp(U_i)+(\varCompGVA)^{-1}}\\
&&\qquad
={\De_i-\chi_i \exp(U_i)-(\varCompGVA)^{-1}
(U_i-\lambdaiGVA/2)\over(B_{i0}\z+\chi_i) \exp
(U_i)+(\varCompGVA)^{-1}}
+O\sfive+O\seight,
\end{eqnarray*}
which implies that
\begin{eqnarray*}
\delta_i&=&{\De_i-\chi_i \exp(U_i)\over(B_{i0}\z+\chi_i) \exp
(U_i)}+O\sfour\\
&=&\{n \exp(\be_0\z)
\phi(\be_1\z)\}\mo\{\De_i \exp(-U_i)-\chi_i\}
+O\sfour\\
&=&\{n \exp(\be_0\z)
\phi(\be_1\z)\}\mo\De_i \exp(-U_i)
-(\betazGVA-\be_0\z)
-(\betaoGVA-\be_1\z) \ga(\be_1\z)+O\sfour.
\end{eqnarray*}
Here we have defined $O\sfour$ is as in Table \ref{tabOnotation}
and have used the fact that
$n\mo B_{i0}\z=\exp(\be_0\z) \phi(\be_1\z)+O\sthree$ and
\[
n\mo B_{i1}\z=\exp(\be_0\z) \phi'(\be_1\z
)+O\sthree
=\exp(\be_0\z) \phi(\be_1\z) \ga(\be
_1\z)+O\sthree.
\]
Therefore,
%
%
\begin{eqnarray}\label{eqTwoPtTwentyThree}
\muiGVA&=&U_i-\tfrac{1}{2}\lambdaiGVA+\delta_i\nonumber\\
&=&U_i+\{n \exp(\be_0\z)
\phi(\be_1\z)\}\mo\De_i \exp(-U_i)\nonumber\\[-8pt]\\[-8pt]
&&{} -(\betazGVA-\be_0\z)
-(\betaoGVA-\be_1\z) \ga(\be_1\z)
+O\sfour\nonumber\\
&=&U_i-\bU+\{n \exp(\be_0\z)
\phi(\be_1\z)\}\mo\De_i \exp(-U_i)+O\sfour,
\nonumber
\end{eqnarray}
where to obtain the second identity we used (\ref{eqTwoPtNineteen})
to place $\lambdaiGVA$ into the
remainder, and to obtain the third identity we used (\ref
{eqTwoPtTwentyTwo}) to show that
$\betazGVA-\be_0\z+(\betaoGVA-\be_1\z) \ga'(\be_1\z)=\bU
+O\sfour$. Here\vspace*{1pt} we have
used the property, deducible from (\ref{eqTwoPtEleven}), (\ref
{eqTwoPtSeventeen})
and (\ref{eqTwoPtNineteen}), that $\ze_i=O\sthree$
and $\bze=O\sfour$.

The next step is to substitute the right-hand side of (\ref
{eqTwoPtTwentyThree})
for $\muiGVA$, and the
right-hand side of (\ref{eqTwoPtNineteen}) for $\lambdaiGVA$, in
(\ref{eqTwoPtEleven}),
and derive an expansion, at (\ref{eqTwoPtTwentySix}) below, of~$\ze_i$.
We obtain
\begin{eqnarray*}
&&[1+\{n \exp(\be_0\z)
\phi(\be_1\z)\}\mo\De_i \exp(-U_i)-\bU]
\oon\sumjon\{\exp(\betaoGVA X\ourij)
-\phi(\betaoGVA)\}\\
&&\qquad=-\phi(\be_1\z) \biggl(\ze_i+\thf\ze_i^2\biggr)
+\oon\sumjon\{Y\ourij\exp(-\be_0\z-U_i)
-\phi(\be_1\z)\}\\
&&\qquad\quad{}
-\exp(-\be_0\z-U_i) (\varCompGVA n)\mo U_i+O\sfive,
\end{eqnarray*}
whence
%
%
\begin{eqnarray}\label{eqTwoPtTwentyFour}\quad
&&\phi(\be_1\z) \biggl(\ze_i+\thf\ze_i^2
\biggr)\nonumber\\
&&\qquad=\oon\sumjon\{Y\ourij\exp(-\be_0\z-U_i)
-\phi(\be_1\z)\}
-\oon\sumjon\{\exp(\betaoGVA X\ourij)
-\phi(\betaoGVA)\}\nonumber\\
&&\qquad\quad{}
-[\{n \exp(\be_0\z)
\phi(\be_1\z)\}\mo\De_i \exp(-U_i)-\bU]
\\
&&\qquad\quad\hspace*{11pt}{} \times\oon\sumjon\{\exp(\betaoGVA X\ourij
)
-\phi(\betaoGVA)\}
\nonumber\\
&&\qquad\quad{}
-\exp(-\be_0\z-U_i) (\varCompGVA n)\mo U_i+O\sfive.
\nonumber
\end{eqnarray}
However, defining
%
%
\begin{equation}\label{eqTwoPtTwentyFive}
D_{ik}(b)=\oon\sumjon\bigl\{X\ourij^k \exp(b X\ourij)
-\phi^{(k)}(b)\bigr\}=O\sthree
\end{equation}
for $k=0,1,2$, and $\De_i=Y\ib-B_{i0}\z\exp(U_i)$, we see that
\begin{eqnarray*}
&&\sumjon\{Y\ourij\exp(-\be_0\z-U_i)
-\phi(\be_1\z)\}
-\sumjon\{\exp(\betaoGVA X\ourij)
-\phi(\betaoGVA)\}\\
&&\qquad=\sumjon\{Y\ourij\exp(-\be_0\z-U_i)
-\phi(\be_1\z)\}\\
&&\qquad\quad{}
-n \{D_{i0}(\be_1\z)
+(\betaoGVA-\be_1\z) D_{i1}(\be_1\z)\}
+O\sthree\\
&&\qquad=\De_i \exp(-\be_0\z-U_i)
-n (\betaoGVA-\be_1\z) D_{i1}(\be_1\z
)+O\sthree,
\end{eqnarray*}
and so, by (\ref{eqTwoPtTwentyFour}),
\begin{eqnarray*}
&&\phi(\be_1\z) \bigl(\ze_i+\tfrac{1}{2}\ze_i^2\bigr)\\
&&\qquad=n\mo\exp(-\be_0\z-U_i) [\De_i
\{1-\phi(\be_1\z)\mo
D_{i0}(\be_1\z)\}-(\varCompGVA)^{-1} U_i]\\
&&\qquad\quad{}
-(\betaoGVA-\be_1\z) D_{i1}(\be_1\z)
+\bU D_{i0}(\be_1\z)+O\sfive.
\end{eqnarray*}
Therefore,
%
%
\begin{eqnarray}\label{eqTwoPtTwentySix}\qquad
\phi(\be_1\z) \ze_i
&=&n\mo\exp(-\be_0\z-U_i) [\De_i
\{1-\phi(\be_1\z)\mo
D_{i0}(\be_1\z)\}-(\varCompGVA)^{-1} U_i
]\nonumber\\
&&{}
-(\betaoGVA-\be_1\z) D_{i1}(\be_1\z)
+\bU D_{i0}(\be_1\z)\\
&&{}
-\tfrac{1}{2}\phi(\be_1\z)\mo
\{n\mo\exp(-\be_0\z-U_i) \De_i\}^2+O\sfive.
\nonumber
\end{eqnarray}

Result (\ref{eqTwoPtTwentySix}), and the fact that $n/m\ra0$ as
$n\rai$, imply that
%
%
\begin{eqnarray}\label{eqTwoPtTwentySeven}
\phi(\beta_1^0) \frac{1}{m} \sum_{i=1}^m U_i \zeta_i &=&
-\frac{1}{mn} \frac{\exp(-\beta_0^0)}{(\sigma^2)^0} \sum_{i=1}^m
U_i^2 \exp(-U_i) \nonumber\\
&&{} -\frac{1}{2m} \phi(\beta_1^0)^{-1}
\sum_{i=1}^m U_i\{ n^{-1} \exp(-\beta_0^0 - U_i)
\Delta_i\}^2 \nonumber\\[-8pt]\\[-8pt]
&&{} + o_p(n^{-1})\nonumber\\
&=& -\frac{1}{n} \exp\biggl\{ {\frac{1}{2}} (\sigma^2)^0 -
\beta_0^0\biggr\} \biggl(1 + {\frac{1}{2}}( \sigma^2)^0
\biggr)
 + o_p(n^{-1}).
\nonumber
\end{eqnarray}
Here we have used the fact that $E\{U_i^2 \exp(-U_i)\}
=\exp(\thf\varCompZero) \varCompZero(1+\varCompZero)$.

\subsection{\texorpdfstring{Initial approximation to
$\varComp-\varCompZero$}{Initial approximation to $sigma^2-(sigma^2)^0$}}

Starting from (\ref{eqTwoPtTwentyOne}), using (\ref{eqTwoPtTwentyTwo})
to substitute for $\betazGVA-\be_0\z$, using
(\ref{eqTwoPtNineteen}) to substitute for $\lambdaiGVA$ and defining
$\bU=m\mo\sumi U_i$ and
$\bze=m\mo\sumi\ze_i$, we obtain
%
%
\begin{eqnarray}\label{eqTwoPtTwentyEight}
\muiGVA&=&U_i+\ze_i-\tfrac{1}{2}\lambdaiGVA
-\log\{\phi(\betaoGVA)/\phi(\be_1\z)\}
-(\betazGVA-\be_0\z)\nonumber\\
&=&U_i+\ze_i-\tfrac{1}{2}\lambdaiGVA-(\betaoGVA-\be_1\z) \ga
(\be_1\z)
-(\betazGVA-\be_0\z)+O\stwo\nonumber\\[-8pt]\\[-8pt]
&=&U_i+\ze_i-
\{2n \phi(\be_1\z) \exp(U_i+\be_0\z
)\}\mo
-(\bU+\bze)\nonumber\\
&&{} +\bigl\{2n \phi(\be_1\z)
\exp\bigl(\be_0\z-\tfrac{1}{2}\varCompZero\bigr)\bigr\}\mo+O\ssix.
\nonumber
\end{eqnarray}
Hence, squaring both sides of (\ref{eqTwoPtTwentyEight}) and
adding,
%
%
\begin{eqnarray}\label{eqTwoPtTwentyNine}\qquad
\oom\sumiom\muiGVA^2
&=&\oom\sumiom(U_i+\ze_i-\bU-\bze)^2\nonumber\\
&&{} -\{mn\phi(\beta_1^0)\exp(\beta_0^0)\}^{-1}
\sum_{i=1}^{m}\exp(-U_i)(U_i+\zeta_i-\bar{U}-\bar{\zeta})\\
&&{} +O\ssix.
\nonumber
\end{eqnarray}

Combining (\ref{eqTwoPtFour}), (\ref{eqTwoPtNineteen}), (\ref
{eqTwoPtTwentySix})
and (\ref{eqTwoPtTwentyNine}) we deduce that
%
%
\begin{eqnarray}\label{eqTwoPtThirty}
\varCompGVA&=&\oom\sumiom(\lambdaiGVA+\muiGVA^2)\nonumber\\
&=&\varCompZero+\oom\sumiom\{(U_i+\ze_i-\bU-\bze
)^2-\varCompZero\}\\
&&{} +\biggl\{n \phi(\be_1\z)
\exp\biggl(\be_0\z-\thf\varCompZero\biggr)\biggr\}\mo
\bigl(1+\varCompZero\bigr)
+O\ssix.
\nonumber
\end{eqnarray}

\subsection{\texorpdfstring{Approximations to $\xi_i$ and $\eta_i$}{Approximations to $xi_i$ and $eta_i$}}

The approximations are given at (\ref{eqTwoPtThirtyOne}) and (\ref
{eqTwoPtThirtyTwo}),
respectively, and are
derived as follows. Note the definition of $D_{ik}(b)$ at (\ref
{eqTwoPtTwentyFive}). In that
notation, observing that $n/m\ra0$ and recalling (\ref
{eqTwoPtFifteen}), it can be deduced
from (\ref{eqTwoPtNine}) and (\ref{eqTwoPtTen}) that, uniformly
in $1\le i\le m$,
%
%
\begin{eqnarray}\label{eqTwoPtThirtyOne}
\quad\xi_i&=&\phi'(\be_1\z)\mo D_{i1}(\be_1\z)
-\tfrac{1}{2}\{\phi'(\be_1\z)\mo
D_{i1}(\be_1\z)\}^2
+O\sfive,
\\
%
\label{eqTwoPtThirtyTwo}
\eta_i&=&\phi'(\be_1\z)\mo[D_{i1}(\be_1\z)\nonumber\\
&&\hphantom{\phi'(\be_1\z)\mo[}{}+(\betaoGVA-\be_1\z)
\{D_{i2}(\be_1\z)
-\phi'(\be_1\z)\mo
\phi''(\be_1\z) D_{i1}(\be_1\z)\}]\\
&&{}
-\tfrac{1}{2}\{\phi'(\be_1\z)\mo
D_{i1}(\be_1\z)\}^2
+O\sfive.
\nonumber
\end{eqnarray}

Result (\ref{eqTwoPtThirtyOne}) is derived by writing (\ref
{eqTwoPtNine}) as
%
%
\begin{equation}\label{eqTwoPtThirtyThree}
\phi'(\be_1\z)\mo D_{i1}(\be_1\z)
=\exp(\xi_i)-1
=\xi_i+\tfrac{1}{2}\xi_i^2+O_p(|\xi_i|^3) ,
\end{equation}
and then inverting the expansion. [The result
$\max_{1\le i\le m} |\xi_i|=o_p(1)$,
in fact $O\sthree$, used in this argument, is readily derived.]
To obtain (\ref{eqTwoPtThirtyTwo}),
note that the analog of (\ref{eqTwoPtThirtyThree}) in that case is
%
%
\begin{equation}\label{eqTwoPtThirtyFour}
\phi'(\betaoGVA)\mo D_{i1}(\betaoGVA)
=\exp(\eta_i)-1=\eta_i+\tfrac{1}{2}\eta_i^2
+O_p(|\eta_i|^3) ,
\end{equation}
and that, uniformly in $1\le i\le m$,
%
%
\begin{eqnarray}\label{eqTwoPtThirtyFive}
&&\phi'(\betaoGVA)\mo D_{i1}(\betaoGVA)\nonumber\\
&&\qquad=\bigl\{\phi'(\be_1\z)
+(\betaoGVA-\be_1\z) \phi''(\be_1\z)
+O\stwo\bigr\}^{ -1}\nonumber\\
&&\qquad\quad{}
\times\bigl\{D_{i1}(\be_1\z)
+(\betaoGVA-\be_1\z) D_{i2}(\be_1\z)+O\sseven
\bigr\}\nonumber\\
&&\qquad=\phi'(\be_1\z)\mo
\{1-(\betaoGVA-\be_1\z)
\phi'(\be_1\z)\mo\phi''(\be_1\z)\}
\\
&&\qquad\quad{}
\times\{D_{i1}(\be_1\z)
+(\betaoGVA-\be_1\z) D_{i2}(\be_1\z)\}
+O\sseven\nonumber\\
&&\qquad=\phi'(\be_1\z)\mo
[D_{i1}(\be_1\z)+(\betaoGVA-\be_1\z)
\{D_{i2}(\be_1\z)-\phi'(\be_1\z)\mo\phi''(\be_1\z)
D_{i1}(\be_1\z)\}]\hspace*{-15pt}\nonumber\\
&&\qquad\quad{}+O\sseven.
\nonumber
\end{eqnarray}

Result (\ref{eqTwoPtThirtyTwo}) follows from (\ref{eqTwoPtThirtyFour})
and (\ref{eqTwoPtThirtyFive}) on inverting the expansion
at (\ref{eqTwoPtThirtyFour}).

\subsection{\texorpdfstring{Another approximation to $\betaoGVA-\be_1\z$, and final
approximations to $\betazGVA-\be_0\z$ and $\varComp-\varCompZero$}{Another approximation to
$beta_1-beta_1^0$, and final
approximations to $beta_0-beta_0^0$ and $sigma^2-(sigma^2)^0$}}
\label{secanother}

Next we use the expansions (\ref{eqTwoPtThirtyOne}), (\ref
{eqTwoPtThirtyTwo}) and (\ref{eqTwoPtTwentySix}) of $\xi_i$, $\eta_i$
and $\ze_i$ to refine the approximations derived in Section 2.3.
The\vspace*{1pt} results are given in (\ref{eqTwoPtFortyTwo}),
(\ref{eqTwoPtFortyThree}) and (\ref{eqTwoPtFortySeven}) in the cases of
$\betazGVA-\be_0\z$, $\betaoGVA-\be_1\z$ and $\varComp-\varCompZero$,
respectively.

It can be deduced from (\ref{eqTwoPtThirtyTwo}) and (\ref
{eqTwoPtTwentySix}) that
%
%
\begin{equation}\label{eqTwoPtThirtySix}
\oom\sumiom\exp(U_i+\eta_i+\ze_i)
=\exp\biggl(\thf\varCompZero\biggr)+O\sthree.
\end{equation}
By (\ref{eqTwoPtThirtyOne}), (\ref{eqTwoPtThirtyTwo}) and (\ref
{eqTwoPtTwentySix}),
%
%
\begin{eqnarray}\label{eqTwoPtThirtySeven}
&&\oom\sumiom\exp(U_i) \{\exp(\xi_i)
-\exp(\eta_i+\ze_i)\}\nonumber\\
&&\qquad=\oom\sumiom\exp(U_i) \biggl[\xi_i
-\eta_i-\ze_i
+\thf\{\xi_i^2-(\eta_i+\ze_i)^2\}\biggr]
+O\sfive\nonumber\\[-8pt]\\[-8pt]
&&\qquad=-\oom\sumiom\exp(U_i) \biggl\{\ze_i
+\thf(2 \eta_i \ze_i+\ze_i^2)\biggr\}
+O\sfive\nonumber\\
&&\qquad\quad{} +O_p\bigl(|\betaoGVA-\be_1\z| n^{\ep
-(1/2)}\bigr).
\nonumber
\end{eqnarray}

Defining $O\snine$ as at Table \ref{tabOnotation}
we deduce from (\ref{eqTwoPtTwentySix}) that
%
%
\begin{eqnarray}\label{eqTwoPtThirtyEight}
\oom\sumiom\exp(U_i) \ze_i
&=&-\thf\phi(\be_1\z)\mt
{1\over mn^2} \sumiom\exp(-2 \be_0\z-U_i) \De_i^2\nonumber\\
&&{}
+O_p\{(mn)\mhf\}+O\sfive\\
&=&-(2n)\mo\phi(\be_1\z)\mo
\exp(-\be_0\z)+O\snine,
\nonumber
\end{eqnarray}
where we have used the fact that $n/m\ra0$ and, since $Y\ib$, conditional
on $\cF_i$, has a Poisson distribution with mean $B_{i0}\z\exp
(U_i)$, then
\begin{eqnarray*}
E\{\exp(-U_i) \De_i^2\}
&=&E[\exp(-U_i) \{Y\ib-E(Y\ib\mi\cF_i)\}^2]\\
&=&E\{\exp(-U_i) \var(Y\ib\mi\cF_i)\}\\
&=&E\{\exp(-U_i) B_{i0}\z\exp(U_i)\}\\
&=&E(B_{i0}\z)=n \exp(\be_0\z) \phi(\be
_1\z) .
\end{eqnarray*}
Similarly,
\begin{eqnarray}\label{eqTwoPtThirtyNine}
\oom\sumiom\exp(U_i) \ze_i^2
&=&\phi(\be_1\z)\mt
{1\over mn^2} \sumiom\exp(-2 \be_0\z-U_i) \De
_i^2+O\snine\nonumber\\[-8pt]\\[-8pt]
&=&n\mo\phi(\be_1\z)\mo
\exp(-\be_0\z)+O\snine.\nonumber
\end{eqnarray}
Moreover, since by (\ref{eqTwoPtThirtyTwo}) and (\ref{eqTwoPtTwentySix}),
\[
\eta_i=\phi'(\be_1\z)\mo D_{i1}(\be_1\z
)+O\sfour,\qquad
\ze_i=\phi(\be_1\z)\mo n\mo\exp(-\be_0\z
-U_i) \De_i
+O\sfour,
\]
and for $k\geq0$,
\[
E\{\exp(U_i) D_{ik}(\be_1\z) \exp(-U_i)
\De_i\}
=E\{D_{ik}(\be_1\z) E(\De_i\mi\cF_i)\}
=0 ,
\]
then
%
%
\begin{equation}\label{eqTwoPtForty}
\oom\sumiom\exp(U_i) \eta_i \ze_i
=O\sfive.
\end{equation}
Together, (\ref{eqTwoPtThirtySeven}),
(\ref{eqTwoPtThirtyEight}), (\ref{eqTwoPtThirtyNine}) and (\ref
{eqTwoPtForty})
imply that
%
%
\begin{eqnarray}\label{eqTwoPtFortyOne}
&&\oom\sumiom\exp(U_i) \{\exp(\xi_i)
-\exp(\eta_i+\ze_i)\}\nonumber\\
&&\qquad=(2n)\mo\phi(\be_1\z)\mo
\exp(-\be_0\z)
-(2n)\mo\phi(\be_1\z)\mo
\exp(-\be_0\z)\nonumber\\[-8pt]\\[-8pt]
&&\qquad\quad{}
+O\snine+O_p\bigl(|\betaoGVA-\be_1\z| n^{\ep-(1/2)}
\bigr)\nonumber\\
&&\qquad=O\snine+O_p\bigl(|\betaoGVA-\be_1\z|
n^{\ep-(1/2)}\bigr) .
\nonumber
\end{eqnarray}

Combining (\ref{eqTwoPtTwenty}), (\ref{eqTwoPtThirtySix}) and (\ref
{eqTwoPtFortyOne}),
and noting that $\De=O_p\{(mn)\mhf\}$
and $n/m\ra0$, we deduce that
%
%
\begin{equation}\label{eqTwoPtFortyTwo}
\betaoGVA-\be_1\z=O\snine.
\end{equation}

Together, (\ref{eqTwoPtTwentyTwo}) and (\ref{eqTwoPtFortyTwo})
imply that
%
%
\begin{equation}\label{eqTwoPtFortyThree}
\betazGVA-\be_0\z=\bU+\bze-c_0 n\mo+o_p(m\mhf+n\mo) ,
\end{equation}
where
\[
c_0=\bigl\{2 \phi(\be_1\z)
\exp\bigl(\be_0\z-\tfrac{1}{2}\varCompZero\bigr)\bigr\}\mo.
\]
Result (\ref{eqPGHa}) of Theorem \ref{thmmainThm} is a direct consequence
of (\ref{eqTwoPtFortyThree}) and the property
%
%
\begin{eqnarray}\label{eqTwoPtFortyFour}
\bar{\zeta} &=& -\frac{1}{m} \sum_{i=1}^m U_i
\{n(\sigma^2)^0 \exp(U_i+\beta_0^0)\phi(\betaoZero)\}
^{-1} \nonumber\\
&&{} - {\frac{1}{2}} \phi(\beta_1^0)^{-2}
E\{n^{-1}
\exp(-\beta_0^0 -U_i) \Delta_i \}^2 + o_p(n^{-1}) \\
&=& c_0 n^{-1} + o_p(n^{-1}).
\nonumber
\end{eqnarray}
Results (\ref{eqTwoPtTwentySix}) and (\ref{eqTwoPtFortyTwo}), and
the property
\[
E\{\exp(-2 U_i) \De_i^2\}
=E\{B_{i0}\z\exp(-U_i)\}
=n \exp\bigl(\be_0\z+\tfrac{1}{2}\varCompZero\bigr) \phi(\be
_1\z) ,
\]
imply that
%
%
\begin{eqnarray}\label{eqTwoPtFortyFive}
\oom\sumiom\ze_i^2
&=&\phi(\be_1\z)\mt
{1\over mn^2} \sumiom\exp(-2 \be_0\z-2 U_i) \De
_i^2+o_p(1)\nonumber\\
&=&n\mo\phi(\be_1\z)\mo
\exp\biggl\{\thf\varCompZero-\be_0\z\biggr\}+o_p(n\mo
)\\
&=&2 c_0 n\mo+o_p(n\mo) .
\nonumber
\end{eqnarray}
By (\ref{eqTwoPtTwentySeven}),
%
%
\begin{eqnarray}\label{eqTwoPtFortySix}
\oom\sumiom U_i \ze_i
&=&-\oon\phi(\be_1\z)\mo\exp\biggl(\thf\varCompZero
-\be_0\z\biggr)
\biggl(1+\smhalf\varCompZero\biggr)\nonumber\\[-8pt]\\[-8pt]
&&{}+o_p(n\mo) .\nonumber
\end{eqnarray}
Together, (\ref{eqTwoPtFortyFour})--(\ref{eqTwoPtFortySix}) give
\begin{eqnarray*}
&& \frac{1}{m} \sum_{i=1}^m \{(U_i + \zeta_i - \bar{U}
- \bar{\zeta})^2 - (\sigma^2)^0 \} \\
&&\qquad = \frac{1}{m} \sum_{i=1}^m \bigl(U_i^2 - (\sigma^2)^0\bigr)
+ \frac{1}{m} \sum_{i=1}^m \zeta_i^2 - \bar{\zeta}^2\\
&&\qquad\quad{} + \frac{2}{m}
\sum_{i=1}^m U_i \zeta_i - 2\bar{U} \bar{\zeta} + O_p(m^{-1})\\
&&\qquad = \frac{1}{m} \sum_{i=1}^m \bigl(U_i^2 - (\sigma^2)^0\bigr)+2n^{-1}c_0
-2n^{-1}c_0 \bigl(2 + (\sigma^2)^0\bigr)\\
&&\qquad\quad{} + o_p(m^{-1/2}+n^{-1})\\
&&\qquad = \frac{1}{m} \sum_{i=1}^m \bigl(U_i^2 - (\sigma^2)^0\bigr)
-2n^{-1}c_0\bigl(1+\varCompZero\bigr)\\
&&\qquad\quad{}+ o_p(m^{-1/2}+n^{-1}).
\end{eqnarray*}
Hence, by (\ref{eqTwoPtThirty}),
%
%
\begin{equation}\label{eqTwoPtFortySeven}
\varCompGVA-\varCompZero=\oom\sumiom\bigl(U_i^2-\varCompZero\bigr)
+o_p(m\mhf+n\mo) .
\end{equation}
Result (\ref{eqPGHc}) of Theorem \ref{thmmainThm} is
a direct consequence of (\ref{eqTwoPtFortySeven}).

\subsection{\texorpdfstring{Final approximation to
$\betaoGVA-\be_1\z$}{Final approximation to $beta_1-beta_1^0$}}

Our first step\vspace*{1pt} is to sharpen the expansion of (\ref{eqTwoPtSix})
at (\ref{eqTwoPtSixteen});
see (\ref{eqTwoPtFiftyTwo}), which leads to (\ref
{eqTwoPtFiftySeven}), the principal
analog of (\ref{eqTwoPtSixteen}).

Recall that
%
%
\begin{eqnarray}\label{eqTwoPtFortyNine}
\De_i&=&Y\ib-\exp(\be_0\z+U_i) \sumjon\exp(\be_1\z
X\ourij)\nonumber\\[-8pt]\\[-8pt]
&=&Y\ib-\exp(U_i) B_{i0}\z.\nonumber
\end{eqnarray}
Also, in view of (\ref{eqTwoPtFortyTwo}) and (\ref{eqTwoPtFortyThree}),
\begin{eqnarray*}
B_i&=&\exp(\betazGVA) \sumjon\exp(\betaoGVA X\ourij)\\
&=&\exp(\be_0\z) \biggl\{1+(\betazGVA-\be_0\z)
+\thf(\betazGVA-\be_0\z)^2
+\osx(\betazGVA-\be_0\z)^3\biggr\}\\
&&{} \times
\sumjon\biggl\{1+(\betaoGVA-\be_1\z) X\ourij
+\thf(\betaoGVA-\be_1\z)^2 X\ourij^2\biggr\}\\
&&\hspace*{26pt}{} \times
\exp(\be_1\z X\ourij)
+O_p(m\mt n+m^{-3/2} n\mhf+m\mo+n^{\ep-3})\\
&=&\exp(\be_0\z)
\sumjon\biggl\{1+(\betazGVA-\be_0\z)
+\thf(\betazGVA-\be_0\z)^2\\
&&\hphantom{\exp(\be_0\z)
\sumjon\biggl\{}{}
+\osx(\betazGVA-\be_0\z)^3+(\betaoGVA-\be_1\z) X\ourij
\\
&&\hphantom{\exp(\be_0\z)
\sumjon\biggl\{}{}{}
+\thf(\betaoGVA-\be_1\z)^2 X\ourij^2
+(\betazGVA-\be_0\z) (\betaoGVA-\be_1\z)
X\ourij\biggr\}
\exp(\be_1\z X\ourij)\\
&&{}
+O_p\bigl(m\mhf n^\ep+n^{\ep-(5/2)}\bigr)\\
&=&\biggl\{1+(\betazGVA-\be_0\z)+\thf(\betazGVA
-\be_0\z)^2
+\osx(\betazGVA-\be_0\z)^3\biggr\} B_{i0}\z\\
&&{}
+\{1+(\betazGVA-\be_0\z)\}
(\betaoGVA-\be_1\z) B_{i1}\z
+\thf(\betaoGVA-\be_1\z)^2 B_{i2}\z+O\sten,
\end{eqnarray*}
where $O\sten$ is defined in Table \ref{tabOnotation}. Hence,
recalling that
$\delta_i=\muiGVA+\thf\lambdaiGVA-U_i$, we see that, for each
$\ep>0$, we have,
uniformly in $1\le i\le n$,
%
%
\begin{eqnarray}\label{eqTwoPtFifty}
&&
Y\ib-B_i \exp(\delta_i+U_i)\nonumber\\
&&\qquad=Y\ib-B_{i0}\z\exp(\delta_i+U_i)\nonumber\\
&&\qquad\quad{} -\bigl[\bigl\{(\betazGVA-\be_0\z)+\tfrac{1}{2}
(\betazGVA-\be_0\z)^2
+\tfrac{1}{6}(\betazGVA-\be_0\z)^3\bigr\} B_{i0}\z\\
&&\qquad\quad\hspace*{16pt}{} +\{1+(\betazGVA-\be_0\z)\}
(\betaoGVA-\be_1\z) B_{i1}\z
+\tfrac{1}{2}(\betaoGVA-\be_1\z)^2 B_{i2}\z\bigr]\nonumber\\
&&\qquad\quad\hspace*{11pt}{} \times\exp(\delta_i+U_i)+O\sten.
\nonumber
\end{eqnarray}
Combining (\ref{eqTwoPtFortyNine}) and (\ref{eqTwoPtFifty}) we
obtain
\begin{eqnarray*}
&&Y\ib-B_i \exp(\delta_i+U_i)\\
&&\qquad=\De_i-\exp(U_i)
\bigl(\{\exp(\delta_i)-1\}
B_{i0}\z\\
&&\qquad\quad\hspace*{63.6pt}{}
+\bigl[\bigl\{(\betazGVA-\be_0\z)+\tfrac{1}{2}
(\betazGVA-\be_0\z)^2\\
&&\qquad\quad\hspace*{130pt}{}+\tfrac{1}{6}(\betazGVA-\be_0\z)^3\bigr\} B_{i0}\z\\
&&\qquad\quad\hspace*{80pt}{}
+\{1+(\betazGVA-\be_0\z)\}
(\betaoGVA-\be_1\z) B_{i1}\z\\
&&\qquad\quad\hspace*{140pt}{} +\tfrac{1}{2}(\betaoGVA-\be_1\z)^2
B_{i2}\z\bigr] \exp(\delta_i)\bigr)\\
&&\qquad\quad{} +O\sten.
\end{eqnarray*}
Therefore, (\ref{eqTwoPtSix}) implies that
\begin{eqnarray*}
(\varCompGVA)^{-1} \muiGVA
&=&\De_i-\exp(U_i)
\bigl(\{\exp(\delta_i)-1\}
B_{i0}\z\\
&&\hspace*{63pt}{}
+\bigl[\bigl\{(\betazGVA-\be_0\z)+\tfrac{1}{2}(\betazGVA
-\be_0\z)^2\\
&&\hspace*{129pt}{}
+\tfrac{1}{6}(\betazGVA-\be_0\z)^3\bigr\} B_{i0}\z\\
&&\hspace*{79pt}{}+\{1+(\betazGVA-\be_0\z)\}
(\betaoGVA-\be_1\z) B_{i1}\z\\
&&\hspace*{139pt}{} +\tfrac{1}{2}(\betaoGVA-\be_1\z)^2 B_{i2}\z\bigr]
\exp(\delta_i)\bigr)\\
&&{}+O\sten,
\end{eqnarray*}
or equivalently,
%
%
\begin{eqnarray}\label{eqTwoPtFiftyOne}
&&\exp(U_i)
\bigl(\{\exp(\delta_i)-1\} B_{i0}\z\nonumber\\
&&\qquad\hspace*{16.5pt}{} +\exp(\delta_i) \bigl[\bigl\{(\betazGVA-\be_0\z
)
+\tfrac{1}{2}(\betazGVA-\be_0\z)^2
+\tfrac{1}{6}(\betazGVA-\be_0\z)^3\bigr\} B_{i0}\z\nonumber\\[-8pt]\\[-8pt]
&&\qquad\hspace*{64.5pt}{} +\{1+(\betazGVA-\be_0\z)\}
(\betaoGVA-\be_1\z) B_{i1}\z
+\tfrac{1}{2}(\betaoGVA-\be_1\z)^2 B_{i2}\z\bigr]
\bigr)\nonumber\\
&&\qquad{} +(\varCompGVA)^{-1} \bigl(\delta_i+U_i-\tfrac{1}{2}\lambdaiGVA
\bigr)
=\De_i+O\sten.\nonumber
\end{eqnarray}

Substituting the far right-hand side of (\ref{eqTwoPtNineteen})
for $\lambdaiGVA$ in (\ref{eqTwoPtFiftyOne}) we deduce that
%
%
\begin{eqnarray}\label{eqTwoPtFiftyTwo}\qquad
&&\exp(\delta_i)-1
+\exp(\delta_i) \bigl\{(\betazGVA-\be_0\z)+\tfrac{1}{2}
(\betazGVA-\be_0\z)^2
+(\betaoGVA-\be_1\z) (B_{i1}\z/B_{i0}\z
)\bigr\}\nonumber\\
&&\quad{}
+\{\varCompGVA B_{i0}\z\exp(U_i)\}\mo
(\delta_i+U_i)\\
&&\qquad=\{B_{i0}\z\exp(U_i)\}\mo\De_i
+O\seleven,
\nonumber
\end{eqnarray}
where $O\seleven$ is as defined in Table \ref{tabOnotation}.
Result (\ref{eqTwoPtFiftyTwo})
implies that
%
%
\begin{equation}\label{eqTwoPtFiftyThree}
\delta_i+\tfrac{1}{2}\delta_i^2 G_{i2}+\tfrac{1}{6}\delta_i^3 G_{i3}
=G_i+O\seleven,
\end{equation}
where, putting
%
%
\begin{eqnarray}\label{eqTwoPtFiftyFour}
G_{i1}
&=&1+(\betazGVA-\be_0\z)+\tfrac{1}{2}(\betazGVA-\be_0\z
)^2
+(\betaoGVA-\be_1\z) (B_{i1}/B_{i0}\z)\nonumber\\[-8pt]\\[-8pt]
&&{}+\{\varCompGVA B_{i0}\z\exp(U_i)\}\mo,\nonumber
\end{eqnarray}
we define $G_i$, $G_{i2}$ and $G_{i3}$ by $G_{i3} G_{i1}=1$,
%
%
\begin{eqnarray}\label{eqTwoPtFiftyFive}
G_{i2} G_{i1}
&=&1+(\betazGVA-\be_0\z)
+(\betaoGVA-\be_1\z) (B_{i1}/B_{i0}\z) ,
\\
%
\label{eqTwoPtFiftySix}
G_i G_{i1}&=&\{B_{i0}\z\exp(U_i)\}\mo\De_i
-\{\varCompGVA B_{i0}\z\exp(U_i)\}\mo U_i\nonumber\\[-8pt]\\[-8pt]
&&{}
-\bigl\{(\betazGVA-\be_0\z)+\tfrac{1}{2}(\betazGVA-\be
_0\z)^2
+(\betaoGVA-\be_1\z) (B_{i1}\z/B_{i0}\z
)\bigr\} .
\nonumber
\end{eqnarray}

Solving (\ref{eqTwoPtFiftyThree}) for $\delta_i$ we deduce that,
for each $\ep>0$,
%
%
\begin{equation}\label{eqTwoPtFiftySeven}
\delta_i=G_i-\tfrac{1}{2} G_{i2} G_i^2
-\bigl(\tfrac{1}{6} G_{i3}-\tfrac{1}{2} G_{i2}^2\bigr) G_i^3
+O\seleven,
\end{equation}
uniformly in $1\le i\le n$. Now, $G_{i1}$, $G_{i2}$ and $G_{i3}$
each equal $1+O_p(m\mhf+n^{\ep-1})$. Therefore, $\osx G_{i3}-\thf G_{i2}^2
=-\otd+O_p(m\mhf+n^{\ep-1})$.
Using (\ref{eqTwoPtFiftyFour}), (\ref{eqTwoPtFiftyFive}) and (\ref
{eqTwoPtFiftySix})
we deduce that
\[
G_{i2}=1-\{\varCompGVA B_{i0}\z\exp(U_i)\}\mo
+O_p(m\mo+n^{\ep-2}) ,\qquad
G_i=H_i+O\seleven,
\]
where
\begin{eqnarray}\label{eqTwoPtFiftyEight}\quad
H_i&=&\bigl[\{B_{i0}\z\exp(U_i)\}\mo\De_i
-\{\varCompGVA B_{i0}\z\exp(U_i)\}\mo U_i\nonumber\\
&&\hspace*{2pt}{}
-\bigl\{(\betazGVA-\be_0\z)+\tfrac{1}{2}(\betazGVA-\be
_0\z)^2
+(\betaoGVA-\be_1\z) (B_{i1}\z/B_{i0}\z
)\bigr\}\bigr]\\
&&{} \times
[1-(\betazGVA-\be_0\z)
-(\betaoGVA-\be_1\z) (B_{i1}/B_{i0}\z)
-\{\varCompGVA B_{i0}\z\exp(U_i)\}\mo] .\nonumber
\end{eqnarray}
Note too that $G_{i2} H_i^2=H_i^2+O_p(m\mhf n^{\ep-1}+n^{\ep-2})$.
Combining the results from (\ref{eqTwoPtFiftySeven}) down we see that
%
%
\begin{equation}\label{eqTwoPtFiftyNine}
\delta_i=H_i-\tfrac{1}{2} H_i^2+\tfrac{1}{3} H_i^3
+O\seleven.
\end{equation}
Note that, as $a\ra0$,
$\exp(a-\thf a^2+\otd a^3)-1=a+O(a^4)$ as $a\ra0$.
This property, (\ref{eqTwoPtFiftyNine}) and the fact
that $H_i^4=O_p(n^{\ep-2})$
imply that
%
%
\begin{equation}\label{eqTwoPtSixty}
\exp(\delta_i)-1=H_i+O\seleven.
\end{equation}

The formula immediately preceding (\ref{eqTwoPtTwenty}) is equivalent to
%
%
\begin{eqnarray}\label{eqTwoPtSixtyOne}
&&\{1+O_p(m^{-1/2}+n\mo)\}
\ga'(\be_1\z) (\betaoGVA-\be_1\z)
\oom\sumiom\exp(U_i+\eta_i+\ze_i)\nonumber\\
&&\qquad=\De\exp(-\be_0\z) \phi(\be_1\z
)\mo\\
&&\qquad\quad{} +\ga(\be_1\z)
\oom\sumiom\{\exp(\xi_i)
-\exp(\eta_i+\ze_i)\} \exp(U_i) .
\nonumber
\end{eqnarray}
Since $\eta_i$ and $\ze_i$ both equal $O\sthree$ [see\vspace*{1pt} (\ref
{eqTwoPtTwentySix})
and (\ref{eqTwoPtThirtyTwo})],
and $m\mo\sum_{i=1}^m\exp(U_i)=E\{\exp(U_1)\}+o_p(1)
=\exp\{\varCompZero/2\}+o_p(1)$, then
(\ref{eqTwoPtSixtyOne}) implies that
%
%
\begin{eqnarray}\label{eqTwoPtSixtyTwo}
&&
\{1+o_p(1)\} \ga'(\be_1\z) (\betaoGVA-\be_1\z
)
\exp\{\varCompZero/2\}
\nonumber\\
&&\qquad=\De\exp(-\be_0\z) \phi(\be_1\z
)\mo\\
&&\qquad\quad{} +\ga(\be_1\z)
\oom\sumiom\{\exp(\xi_i)
-\exp(\eta_i+\ze_i)\} \exp(U_i) .\nonumber
\end{eqnarray}
Formulae (\ref{eqTwoPtNine}) and (\ref{eqTwoPtTen}) are together
equivalent to
%
%
\begin{eqnarray}
\label{eqTwoPtSixtyThree}
\phi'(\be_1\z)
\{\exp(\xi_i)-1\}
&=&\oon\sumjon\{X\ourij\exp(\be_1\z X\ourij)
-\phi'(\be_1\z)\} ,
\\
\label{eqTwoPtSixtyFour}
\phi'(\betaoGVA) \{\exp(\eta_i)-1\}
&=&\oon\sumjon\{X\ourij\exp(\betaoGVA X\ourij)
-\phi'(\betaoGVA)\} .
\end{eqnarray}
Result (\ref{eqTwoPtSixtyFour}) implies that, for each $\ep>0$,
\begin{eqnarray*}
&&\{\phi'(\be_1\z)
+O_p(|\betaoGVA-\be_1\z|)\} \{\exp(\eta
_i)-1\}\\
&&\qquad=\oon\sumjon\{X\ourij\exp(\be_1\z X\ourij
)
-\phi'(\be_1\z)\}
+O_p\bigl(|\betaoGVA-\be_1\z| n^{\ep-(1/2)}\bigr) ,
\end{eqnarray*}
uniformly in $1\le i\le n$. Therefore, since $\eta_i=O\sthree$
[see (\ref{eqTwoPtThirtyTwo})], then
\begin{eqnarray*}
\phi'(\be_1\z) \{\exp(\eta_i)-1\}
&=&\oon\sumjon\{X\ourij\exp(\be_1\z X\ourij)
-\phi'(\be_1\z)\}\\
&&{}+O_p\bigl(|\betaoGVA-\be_1\z| n^{\ep-(1/2)}\bigr),
\end{eqnarray*}
which in company with (\ref{eqTwoPtSixtyFour}) implies that
\[
\phi'(\be_1\z) \{\exp(\eta_i)-\exp(\xi_i)\}
=O_p\bigl(|\betaoGVA-\be_1\z| n^{\ep-(1/2)}\bigr) ,
\]
uniformly in $1\le i\le n$. Hence, since $\eta_i=O\sthree$ and
$\ze_i=O\sthree$ [see (\ref{eqTwoPtTwentySix}) and (\ref
{eqTwoPtThirtyTwo})],
%
%
\begin{eqnarray}\label{eqTwoPtSixtyFive}
\exp(\xi_i)-\exp(\eta_i+\ze_i)
&=&\{\exp(\xi_i)-\exp(\eta_i)\} \exp(\ze_i)\nonumber\\
&&{} +\exp(\xi_i) \{1-\exp(\ze_i)\}\nonumber\\[-8pt]\\[-8pt]
&=&\exp(\xi_i) \{1-\exp(\ze_i)\}\nonumber\\
&&{} +O_p\bigl(|\betaoGVA-\be_1\z| n^{\ep-(1/2)}
\bigr) ,
\nonumber
\end{eqnarray}
uniformly in $i$. Combining (\ref{eqTwoPtSixtyTwo}) and (\ref
{eqTwoPtSixtyFive})
we deduce that
%
%
\begin{eqnarray}\label{eqTwoPtSixtySix}
&&\{1+o_p(1)\} \ga'(\be_1\z) (\betaoGVA-\be_1\z
)
\exp\biggl\{\smhalf\varCompZero\biggr\}
\nonumber\\
&&\qquad=\De\exp(-\be_0\z) \phi(\be_1\z
)\mo\\
&&\qquad\quad{} +\ga(\be_1\z)
\oom\sumiom\exp(\xi_i+U_i) \{1-\exp(\ze_i)\} .
\nonumber
\end{eqnarray}

Next we return to (\ref{eqTwoPtEleven}), which we write equivalently as
%
%
\begin{eqnarray}\label{eqTwoPtSixtySeven}
\phi(\be_1\z) \{1-\exp(\ze_i)\}
&=&\exp(\betazGVA-\be_0\z+\delta_i)
\oon\sumjon\{\exp(\betaoGVA X\ourij)-\phi(\betaoGVA)\}
\nonumber\\
&&{} -\oon\sumjon\{Y\ourij\exp(-\be
_0\z-U_i)
-\phi(\be_1\z)\}\\
&&{} +(\varCompGVA n)\mo\muiGVA\exp(-\be
_0\z-U_i).
\nonumber
\end{eqnarray}
So that we might replace $\betaoGVA$ by $\be_1\z$ on the right-hand
side of
(\ref{eqTwoPtSixtySeven}), we observe that
%
%
\begin{eqnarray}\label{eqTwoPtSixtyEight}
\oon\sumjon\{\exp(\betaoGVA X\ourij)-\phi(\betaoGVA)\}
&=&\oon\sumjon\{\exp(\be_1\z X\ourij)
-\phi(\be_1\z)\}
\nonumber\\[-8pt]\\[-8pt]
&&{}
+O_p\bigl(|\betaoGVA-\be_1\z| n^{\ep-(1/2)}\bigr) .
\nonumber
\end{eqnarray}
Combining (\ref{eqTwoPtSixtySix})--(\ref{eqTwoPtSixtyEight}) we obtain
%
%
\begin{eqnarray}\label{eqTwoPtSixtyNine}
&&\{1+o_p(1)\} \ga'(\be_1\z) (\betaoGVA-\be_1\z
)
\exp\biggl\{\smhalf\varCompZero\biggr\}
\nonumber\\
&&\qquad=\De\exp(-\be_0\z) \phi(\be_1\z)\mo
\nonumber\\
&&\qquad\quad{}
+{\phi'(\be_1\z)\over\phi(\be_1\z)^2}
\oom\sumiom\exp(\xi_i+U_i) \nonumber\\[-8pt]\\[-8pt]
&&\hspace*{105.7pt}{}\times\Biggl[\exp(\betazGVA-\be_0\z
+\delta_i)
\oon\sumjon\{\exp(\be_1\z X\ourij)
-\phi(\be_1\z)\}\nonumber\\
&&\qquad\quad\hspace*{90pt}{} -\oon\sumjon\{Y\ourij\exp(-\be_0\z-U_i)
-\phi(\be_1\z)\}\nonumber\\
&&\hspace*{201.6pt}{}+(\varCompGVA n)\mo\muiGVA\exp(-\be_0\z-U_i)\Biggr] .
\nonumber
\end{eqnarray}
(Recall that $\ga=\phi'\phi\mo$, and so $\ga/\phi=\phi'\phi\mt$.)

Since $\exp(\xi_i)-1=D_{i1}(\be_1\z) \phi'(\be_1\z)\mo$
[see (\ref{eqTwoPtNine})]
and $\betazGVA-\be_0\z=O_p(m\mhf+n\mo)$ [see (\ref
{eqTwoPtFortyThree})], then
%
%
\begin{eqnarray}\label{eqTwoPtSeventy}
&&\oom\sumiom\exp(\xi_i+U_i) \exp(\betazGVA-\be_0\z
+\delta_i)
\oon\sumjon\{\exp(\be_1\z X\ourij)
-\phi(\be_1\z)\}\nonumber\\
&&\qquad=\biggl\{1+(\betazGVA-\be_0\z)
+\thf(\betazGVA-\be_0\z)^2\biggr\}
\oom\sumiom\exp(\xi_i+\delta_i+U_i) D_{i0}(\be_1\z
)\nonumber\\
&&\qquad\quad{}
+O_p(m^{-3/2}+n^{-3})\\
&&\qquad=\biggl\{1+(\betazGVA-\be_0\z)
+\thf(\betazGVA-\be_0\z)^2\biggr\}\nonumber\\
&&\qquad\quad{} \times
\oom\sumiom\exp(\delta_i+U_i)
\{1+D_{i1}(\be_1\z) \phi'(\be_1\z)\mo\}
D_{i0}(\be_1\z)\nonumber\\
&&\qquad\quad{}
+O_p(m^{-3/2}+n^{-3}).
\nonumber
\end{eqnarray}
Likewise,
%
%
\begin{eqnarray}\label{eqTwoPtSeventyOne}
&&\oom\sumiom\exp(\xi_i+U_i)
\oon\sumjon\{Y\ourij\exp(-\be_0\z-U_i)
-\phi(\be_1\z)\}\nonumber\\
&&\qquad={1\over m} \sumiom\exp(U_i)
\{1+D_{i1}(\be_1\z) \phi'(\be_1\z)\mo
\}\\
&&\qquad\quad\hphantom{{1\over m} \sumiom}{}\times\{n\mo\De_i \exp(-\be_0\z-U_i)
+D_{i0}(\be_1\z)\}
\nonumber
\end{eqnarray}
and, since $\sumi\muiGVA=0$ [see (\ref{eqTwoPtSeven})],
%
%
\begin{eqnarray}\label{eqTwoPtSeventyTwo}
&&\oom\sumiom\exp(\xi_i+U_i)
(\varCompGVA n)\mo\muiGVA\exp(-\be_0\z-U_i)\nonumber
\\[-0.5pt]
&&\qquad={1\over\varCompGVA mn}
\sumiom\exp(\xi_i-\be_0\z) \muiGVA\nonumber\\[-0.5pt]
&&\qquad=\exp(-\be_0\z) {1\over\varCompGVA mn}
\sumiom\{1+D_{i1}(\be_1\z) \phi'(\be_1\z
)\mo\}
\muiGVA\\[-0.5pt]
&&\qquad=\exp(-\be_0\z) \phi'(\be_1\z)\mo
{1\over\varCompGVA mn}
\sumiom D_{i1}(\be_1\z) \muiGVA\nonumber\\[-0.5pt]
&&\qquad=O_p(m\mhf n^{-3/2}) .
\nonumber
\end{eqnarray}

Combining (\ref{eqTwoPtSixtyNine})--(\ref{eqTwoPtSeventyTwo}) we
see that
%
%
\begin{eqnarray}\label{eqTwoPtSeventyThree}
&&\{1+o_p(1)\} \ga'(\be_1\z) (\betaoGVA-\be_1\z
)
\exp\biggl\{\smhalf\varCompZero\biggr\}\nonumber\\[-0.5pt]
&&\qquad=\De\exp(-\be_0\z) \phi(\be_1\z
)\mo\nonumber\\[-0.5pt]
&&\qquad\quad{}
+{\phi'(\be_1\z)\over\phi(\be_1\z)^2}
\biggl[\biggl\{1+(\betazGVA-\be_0\z)+\thf(\betazGVA-\be_0\z)^2\biggr\}\nonumber\\[-0.5pt]
&&\qquad\quad\hspace*{52pt}{} \times
\oom\sumiom\exp(\delta_i+U_i)
\{1+D_{i1}(\be_1\z) \phi'(\be_1\z)\mo
\}
D_{i0}(\be_1\z)\\[-0.5pt]
&&\qquad\quad\hspace*{53pt}{}
-\exp(-\be_0\z) {1\over m} \sumiom
\{1+D_{i1}(\be_1\z) \phi'(\be_1\z)\mo
\}\nonumber\\[-0.5pt]
&&\hspace*{130pt}\qquad\quad{} \times
\{n\mo\De_i
+\exp(\be_0\z+U_i) D_{i0}(\be_1\z)\}
\biggr]\nonumber\\[-0.5pt]
&&\qquad\quad{}
+O_p(m\mhf n\mo+n^{-3}) .
\nonumber
\end{eqnarray}

Using the fact that $E(\De_i\mi\cF_i)=0$ and $D_{i1}(\be_1\z
)=O\sthree$
it can be proved that, for all $\ep>0$,
%
%
\begin{eqnarray}\label{eqTwoPtSeventyFour}
&&{1\over mn} \sumiom\exp(-\be_0\z)
\{1+D_{i1}(\be_1\z) \phi'(\be_1\z)\mo
\} \De_i\nonumber\\[-8pt]\\[-8pt]
&&\qquad=\exp(-\be_0\z) {1\over
mn} \sumiom\De_i
+O_p(m\mhf n\mo) .
\nonumber
\end{eqnarray}
Also,
%
%
\begin{eqnarray}\label{eqTwoPtSeventyFive}
\De'&\equiv&\De\exp(-\be_0\z) \phi(\be_1\z
)\mo
-\frac{\exp(-\be_0\z)\phi'(\be_1\z)}{\phi(\be_1\z
)^2} {1\over mn} \sumiom\De_i\nonumber\\[-0.5pt]
&=&\phi(\be_1\z)^{-1}\exp(-\be_0\z) {1\over
mn} \\[-0.5pt]
&&{}\times\sumiom\sumjon
\biggl\{X\ourij-\frac{\phi'(\be_1\z)}{\phi(\be_1\z
)}\biggr\}
\{Y\ourij-\exp(\be_0\z+\be_1\z X\ourij+U_i)
\}.
\nonumber
\end{eqnarray}
Moreover, using (\ref{eqTwoPtFortyThree}) and the fact that
$D_{i0}(\be_1\z)=O\sthree$ and
$E\{D_{i0}(\be_1\z)\mi U_i\}=0$, it can be shown that
%
%
\begin{eqnarray}\label{eqTwoPtSeventySix}
&&(\betazGVA-\be_0\z)
\oom\sumiom\exp(U_i)
\{1+D_{i1}(\be_1\z) \phi'(\be_1\z)\mo
\}
D_{i0}(\be_1\z)\nonumber\\[-0.5pt]
&&\qquad= O_p\bigl\{(m\mhf+n\mo)\cdot
\bigl(m\mhf n^{\ep-(1/2)}\bigr)\bigr\}\\[-0.5pt]
&&\qquad=O_p(m\mhf n^{\ep-1}) .
\nonumber
\end{eqnarray}

Combining (\ref{eqTwoPtSeventyThree})--(\ref{eqTwoPtSeventySix}) we
deduce that
%
%
\begin{eqnarray}\label{eqTwoPtSeventySeven}
&&\{1+o_p(1)\} \ga'(\be_1\z) (\betaoGVA-\be_1\z
)
\exp\biggl\{\smhalf\varCompZero\biggr\}\nonumber\\
&&\qquad=\De'
+{\phi'(\be_1\z)\over\phi(\be_1\z)^2}
\oom\sumiom\exp(U_i) \{\exp(\delta_i)-1\} \nonumber\\[-8pt]\\[-8pt]
&&\qquad\quad\hspace*{84.6pt}{}
\times\{1+D_{i1}(\be_1\z) \phi'(\be_1\z
)\mo\}
D_{i0}(\be_1\z)\nonumber\\
&&\qquad\quad{} +O_p(m\mhf n^{\ep-1}+n^{-3}) .
\nonumber
\end{eqnarray}

Using (\ref{eqTwoPtSixty}) to substitute for $\exp(\delta_i)-1$ in
(\ref{eqTwoPtSeventySeven}), and noting
that $D_{ik}(\be_1\z)=O\sthree$ for $k=0,1$, we deduce from
(\ref{eqTwoPtSeventySeven}) that
%
%
\begin{eqnarray}\label{eqTwoPtSeventyEight}
&&\{1+o_p(1)\} \ga'(\be_1\z) (\betaoGVA-\be_1\z
)
\exp\biggl\{\smhalf\varCompZero\biggr\}\nonumber\\[-8pt]\\[-8pt]
&&\qquad=\De'+{\phi'(\be_1\z)\over\phi(\be
_1\z)^2} \psi(H)
+O_p\bigl(m\mhf n^{\ep-1}+n^{\ep-(5/2)}\bigr) ,
\nonumber
\end{eqnarray}
where $H=(H_1,\ldots,H_m)$, $H_i$ is
as defined at (\ref{eqTwoPtFiftyEight}), and, given a
sequence of random variables $K=(K_1,\ldots,K_m)$, we put
\[
\psi(K)=\oom\sumiom\exp(U_i) K_i
\{1+D_{i1}(\be_1\z) \phi'(\be_1\z)\mo
\}
D_{i0}(\be_1\z) .
\]
Note again that $|D_{i0}(\be_1\z)|=O\sthree$, and the dominant term on
the right-hand side of formula (\ref{eqTwoPtFiftyEight})
for $H_i$ is $\{B_{i0}\z\exp(U_i)\}\mo\De_i$.
Moreover, $|\betazGVA-\be_0\z|=O_p(m\mhf+n\mo)$
[see (\ref{eqTwoPtFortyThree})],
$|\betaoGVA-\be_1\z|=O_p\{(mn)\mhf+n^{\ep-(3/2)}\}$
[see (\ref{eqTwoPtFortyTwo})],
\[
\{\varCompGVA B_{i0}\z\exp(U_i)\}\mo
=\{n \varCompZero\phi(\be_1\z)
\exp(\be_0\z+U_i)\}\mo+O_p\bigl(n^{\ep-(3/2)}
\bigr)
\]
and
\[
B_{i1}\z B_{i0}\z{}\mo
=\phi'(\be_1\z) \phi(\be_1\z)\mo
+O\sthree.
\]
Combining these properties we deduce that (\ref{eqTwoPtSeventyEight})
continues to hold if, on the right-hand side, $\psi(H)$ is replaced by
$\psi(H')$ where $H'=(H_1',\ldots,H_m')$ and
$H_i'=H_i^{(1)}-H_i^{(2)}-H_i^{(3)}$,
with
\begin{eqnarray*}
H_i^{(1)}&=&\{B_{i0}\z\exp(U_i)\}\mo\De_i
[1-(\betazGVA-\be_0\z)
-\{n \varCompZero\phi(\be_1\z) \exp(U_i)\}
\mo] ,\\
H_i^{(2)}&=&\{\varCompGVA B_{i0}\z\exp(U_i)\}\mo U_i
\end{eqnarray*}
and
\[
H_i^{(3)}=(\betazGVA-\be_0\z)+\tfrac{1}{2}
(\betazGVA-\be_0\z)^2
+(\betaoGVA-\be_1\z)
\{\phi'(\be_1\z)/\phi(\be_1\z)\} .
\]
(Note that $H_i^{(3)}$ does not depend on $i$.) It can be proved from the
properties $E(\De_i\mi\cF_i)=0$ and $|D_{i0}(\be_1\z)|=O\sthree$ that,
with $H^{(j)}$ denoting $(H_1^{(j)},\ldots,H_m^{(j)})$, we have
%
%
\begin{equation}\label{eqTwoPtEighty}
\psi\bigl(H^{(1)}\bigr)=O_p(m\mhf n\mo) .
\end{equation}
More simply, since $E(U_i\mi X_{i1},\ldots,X_{in})=0$, then
\begin{eqnarray}\label{eqTwoPtEightyOne}
\psi\bigl(H^{(2)}\bigr)
&=&\oom\sumiom(\varCompGVA B_{i0}\z)\mo U_i
\{1+D_{i1}(\be_1\z) \phi'(\be_1\z)\mo
\}
D_{i0}(\be_1\z)\nonumber\\[-8pt]\\[-8pt]
&=&O_p(m\mhf n^{-3/2}).\nonumber
\end{eqnarray}
Furthermore, writing ${\mathbf1}=(1,\ldots,1)$, an $n$-vector, and
noting that
the properties $E\{D_{ik}(\be_1\z)\mi U_i\}=0$, $\var\{D_{ik}(\be
_1\z)\mi U_i\}=O(n\mo)$
and $E\{\exp(U_i)\}=\exp(\thf\varCompZero)$ imply that
\begin{eqnarray*}
\psi({\mathbf1})&=&\oom\sumiom\exp(U_i)
\{1+D_{i1}(\be_1\z) \phi'(\be_1\z)\mo
\}
D_{i0}(\be_1\z)\\
&=&\phi'(\be_1\z)\mo
\oom\sumiom\exp(U_i) D_{i1}(\be_1\z)
D_{i0}(\be_1\z)
+O_p(m\mhf n\mhf)\\
&=&n\mo\{\phi'(2\be_1\z)
\phi'(\be_1\z)^{-1}-\phi(\be_1\z)\} \exp
\biggl(\thf\varCompZero\biggr)\\
&&{}+O_p(m\mhf n\mhf+n^{-3/2}) ;
\end{eqnarray*}
we obtain
\begin{eqnarray}\label{eqTwoPtEightyTwo}\quad
\psi\bigl(H^{(3)}\bigr)
&=&\biggl[(\betazGVA-\be_0\z)+\thf(\betazGVA-\be
_0\z)^2
+(\betaoGVA-\be_1\z)
\{\phi'(\be_1\z)/\phi(\be_1\z)
\}\biggr]
\psi({\mathbf1})\nonumber\\
&=& \biggl\{\frac{1}{n} [\phi^{\prime}(2\beta^0_1) \phi^{\prime}(\beta_1^0)^{-1}
- \phi(\beta_1^0) ] \exp\biggl( {\frac{1}{2}}
(\sigma^2)^0\biggr) \nonumber\\
&&\hspace*{73pt}{}+ O_p(m^{-1/2} n^{-1/2} + n^{-3/2}) \biggr\} \\
&&{} \times\biggl[ (\widehat{\underline{\beta}}{}_0 - \beta_0^0)
+ \frac{1}{2} (\widehat{\underline{\beta}}{}_0 - \beta_0^0)^2 +
(\widehat{\underline{\beta}}{}_1
- \beta_1^0) \{\phi'(\beta_1^0)/\phi(\beta_1^0) \}\biggr] \nonumber\\
&=& O_p(m^{-1/2} n^{-1}).\nonumber
\end{eqnarray}
To obtain the last line here we used (\ref{eqPGHa})
of Theorem \ref{thmmainThm}, already proved
in Section~\ref{secanother} above.

Combining (\ref{eqTwoPtEighty})--(\ref{eqTwoPtEightyTwo}),
and noting that the function $\psi$ is linear,
so that
\[
\psi(H)=\psi\bigl(H^{(1)}\bigr)-\psi\bigl(H^{(2)}\bigr)-\psi\bigl(H^{(3)}\bigr),
\]
we deduce that
\begin{equation}\label{eqTwoPtEightyThree}
\{1+o_p(1)\} \ga'(\be_1\z) (\betaoGVA-\be_1\z)
\exp\bigl(\tfrac{1}{2}\varCompZero\bigr)
=\De'+o_p\{(mn)\mhf+n\mt\} .\hspace*{-32pt}%
\end{equation}
Furthermore, the random variable $\De'$, defined at (\ref
{eqTwoPtSeventyFive}),
is asymptotically normally distributed with zero mean and variance
\begin{eqnarray*}
&&\frac{\exp(-2 \be_0\z)}{mn}
E\biggl(\biggl\{X_{11}-\frac{\phi'(\be_1\z)}{\phi(\be_1\z)}\biggr\}
^2
E[E\{Y_{11}-E(Y_{11}\mi X_{11},U_1)\}^2
\Bigmi X_{11},U_1]\biggr)\\
&&\qquad=(mn)\mo\exp(-2\be_0\z)
E\biggl[\biggl\{X_{11}-\frac{\phi'(\be_1\z)}{\phi(\be_1\z)}\biggr\}
^2
\exp(\be_0\z+\be_1\z X_{11}+U_1)\biggr]\\
&&\qquad=(mn)\mo\exp\biggl(\thf\varCompZero-\be_0\z\biggr)
E\biggl[\biggl\{X_{11}-\frac{\phi'(\be_1\z)}{\phi(\be_1\z)}\biggr\}
^2
\exp(\be_1\z X_{11})\biggr]\\
&&\qquad=(mn)\mo\ga'(\be_1\z)^2 \exp\{\varCompZero\}
\tau^2 ,
\end{eqnarray*}
where $\tau^2$ is as at (\ref{eqonePtFive}). Result (\ref{eqPGHb})
of the Theorem \ref{thmmainThm} is implied by this property
and (\ref{eqTwoPtEightyThree}).

\section*{Acknowledgments}

The authors are grateful to John Ormerod and Mike Titterington for
their assistance in the preparation of this paper.


%

\printaddresses

\end{document}